\documentclass[12pt]{article}
\usepackage{amsfonts}
\usepackage{mathrsfs}
\usepackage{amsmath}
\usepackage{latexsym}
\usepackage{amssymb}
\usepackage{amsthm}
\usepackage{amscd}

\title{Classification of {\rm PM} Quiver Hopf Algebras }
\author{ \small
Shouchuan Zhang $^{a,~b}$,   Yao-Zhong Zhang $^b$, Hui-xiang Chen $^c$   \ \  \ \\
\small $a$. Department  of Mathematics,
Human University\\  \small   Changsha  410082, \ P.R. China \\
\small $b$. Department of Mathematics, University of Queensland\\
\small Brisbane 4072, Australia\\
\small $c$. Department of Mathematics, Yangzhou University \\
 \small   Yangzhou 225002, P.R. China}
\date{}

\begin{document}
\newtheorem{Proposition}{\quad Proposition}[section]
\newtheorem{Theorem}{\quad Theorem}
\newtheorem{Definition}[Proposition]{\quad Definition}
\newtheorem{Corollary}[Proposition]{\quad Corollary}
\newtheorem{Lemma}[Proposition]{\quad Lemma}
\newtheorem{Example}[Proposition]{\quad Example}
\maketitle \addtocounter{section}{-1}

\numberwithin{equation}{section}

\date{}

\begin {abstract} We describe certain quiver Hopf algebras by parameters. This leads to
the classification of multiple Taft algebras as well as pointed
Yetter-Drinfeld modules and their corresponding Nichols algebras. In
particular, when the ground-field $k$ is the complex field and $G$
is a finite abelian group, we classify quiver Hopf algebras over
$G$, multiple Taft algebras over $G$ and Nichols algebras in
$^{kG}_{kG} {\cal YD}$. We show that the quantum enveloping algebra
of a complex semisimple Lie algebra is a quotient of a semi-path
Hopf algebra.


keywords: Quiver, Hopf algebra, Hopf bimodule
\end {abstract}

\section*{\bf Introduction}

In noncommutative algebras quivers and associated path algebras are
intensively studied because of a remarkable interaction between
homology theory, algebraic geometry, and Lie theory. Ringel's
approach to quantum groups via quivers suggested that there may be
an interesting overlapping also between the theory of quivers and
Hopf algebras. Some years ago Cibils and Rosso   \cite {CR02, CR97}
started to study quivers admitting a (graded) Hopf algebra structure
and there are many papers (e.g. \cite {CHYZ04, OZ04} ) follow these
works. The representation theory of (graded) Hopf algebras has
applications in many branches of physics and mathematics such as the
construction of solutions to the quantum Yang-Baxter equation and
topological invariants (see e.g. \cite{Ma90,RT90}. Quivers have also
been used in gauge theory and string theory (see e.g.
\cite{Zh05,RR04}).

This paper can be viewed an extension of the analysis of Cibils and
Rosso. Let $G$ be a group and $kG$ be the group algebra of $G$ over
a field $k$. It is well-known \cite{CR97} that the $kG$-Hopf
bimodule category $^{kG}_{kG}{\mathcal M}^{kG}_{kG}$ is equivalent
to the direct product category $\prod_{C\in{\mathcal K}(G)}{\mathcal
M}_{kZ_{u(C)}}$, where ${\mathcal K}(G)$ is the set of conjugate
classes in $G$, $u:{\mathcal K}(G)\rightarrow G$ is a map such that
$u(C)\in C$ for any $C\in {\mathcal K}(G)$, $Z_{u(C)}=\{g\in G\mid
gu(C)=u(C)g\}$ and ${\mathcal M}_{kZ_{u(C)}}$ denotes the category
of right $kZ_{u(C)}$-modules (see \cite [Proposition 3.3]{CR97}
\cite [Theorem 4.1]{CR02} and Theorem \ref{1.1}). Thus for any Hopf
quiver $(Q,G,r)$, the $kG$-Hopf bimodule structures on the arrow
comodule $kQ_1^c$ can be derived from the right $kZ_{u(C)}$-module
structures on $^{u(C)}\! (kQ_1^c)^1$ for all $C\in {\mathcal K}(G)$.
If the arrow comodule $kQ_1^c$ admits a $kG$-Hopf bimodule
structure, then there exist six   graded Hopf algebras: co-path Hopf
algebra $kQ^c$,   one-type-co-path  Hopf algebras $kG [kQ_1^c]$,
one-type-path  Hopf algebras $(kG)^* [kQ_1^a]$,  semi-path Hopf
algebra $kQ^s$, semi-co-path Hopf algebra $kQ^{sc}$ and path Hopf
algebra $kQ^a$. We call these Hopf algebras quiver Hopf algebras
over $G$. If the corresponding $kZ_{u(C)}$-modules $^{u(C)}\!
(kQ_1^c)^1$ is pointed (i.e. it is zero or a direct sum of one
dimensional $kZ_{u(C)}$-modules) for all $C\in {\mathcal K}(G)$,
then $kQ_1^c$ is called a {\rm PM} $kG$-Hopf bimodule. The six
graded Hopf algebras derived from the {\rm PM} $kG$-Hopf bimodule
$kQ_1^c$ are called {\rm PM} quiver Hopf algebras. A Yetter-Drinfeld
 $H$-module $V$  is {\it pointed} if $V=0$ or $V$ is a direct
sum of one dimensional {\rm YD} $H$-modules. If $V$ is a pointed
{\rm YD} $H$-module, then the corresponding Nichols algebra
${\mathcal B}(V)$ is called a {\rm PM} {\it Nichols algebra}. For
example, when $G$ is a finite abelian group of exponent $m$ and $k$
contains a primitive $m$-th root of 1 (e.g. \ $k$ is the complex
field ), all quiver Hopf algebras are {\rm PM} ( see Lemma \ref
{1.2}), all {\rm YD } $kG$-modules are
 pointed and all Nichols algebras of all {\rm YD}
$kG$-modules are {\rm PM}
 (see Lemma \ref {2.3}).

The aim of this paper is to provide parametrization of {\rm PM}
quiver Hopf algebras, multiple Taft algebras and {\rm PM} Nichols
algebras. In other words, this paper provides a kind of
classification of these Hopf algebras.
 This paper is organized as follows. In
Section \ref {s1}, we examine the {\rm PM} quiver Hopf algebras by
means of ramification system with characters. In Section \ref {s2},
we describe {\rm PM} Nichols algebras and multiple Taft algebras by
means of element system with characters. In Section \ref {s3}, we
show that the diagram of a quantum weakly commutative multiple Taft
algebra is not only a Nichols algebra but also a quantum linear
space in $^{kG}_{kG}{\cal YD}$; the diagram of a semi-path Hopf
algebra of ${\rm ESC}$   is  a quantum tensor algebra in
$^{kG}_{kG}{\cal YD}$; the quantum enveloping algebra of a complex
semisimple Lie algebra is a quotient of a semi-path Hopf algebra.

\section{\bf Preliminaries}\label{s0}

Throughout this paper, we work over a fixed field $k$. All algebras,
coalgebras, Hopf algebras, and so on, are defined over $k$; dim,
$\otimes$ and Hom stand for ${\rm dim}_k$, $\otimes_k$ and ${\rm
Hom}_k$, respectively. Books \cite {DNR01,Mo93,Sw69} provide the
necessary background for Hopf algebras and book \cite {ARS95}
provides a nice description of the path algebra approach.

Let $\mathbb{Z}$, $\mathbb{Z}^+$ and $\mathbb{N}$ denote sets of all
integers, all positive integers and all non-negative integers,
respectively. For sets $X$ and $Y$, we denote by $|X|$ the cardinal
number of $X$ and by $X^Y$ or $X^{\mid \! Y \! \mid }$ the Cartesian
product $ \Pi _{y \in Y} X_y$ with $X_y = X$ for any $y \in Y$. If
$X$ is finite, then $|X|$ is the number of elements in $X$. If $X =
\oplus _{i\in I} X_{(i)}$ as vector spaces, then we denote by $\iota
_i$ the natural injection from $X_{(i)}$ to $X$ and by $\pi _i$ the
corresponding projection from $X$ to $X_{(i)}$. We will use $\mu$ to
denote the multiplication  of an algebra, $\Delta$ to denote the
comultiplication of a coalgebra, $\alpha ^-$, $\alpha ^+$, $\delta
^-$ and $\delta ^+$ to denote the left module, right module, left
comodule and right comodule structure maps, respectively. The
Sweedler's sigma notations for coalgebras and comodules are $\Delta
(x) = \sum x_{(1)}\otimes x_{(2)}$, $\delta ^- (x)= \sum x_{(-1)}
\otimes x_{(0)}$, $\delta ^+ (x)= \sum x_{(0)} \otimes x_{(1)}$. Let
$G$ be a group. We  denote by $Z(G)$ the center of $G$. Let
$\widehat G$ denote  the set  of characters of all one-dimensional
representations  of $G$. It is clear that $\widehat G$ = $\{ \chi
\mid $  $\chi $ is a group homomorphism from $G$ to the
multiplicative group of all non-zero elements in $k$ \}.

A quiver $Q=(Q_0,Q_1,s,t)$ is an oriented graph, where  $Q_0$ and
$Q_1$ are the sets of vertices and arrows, respectively; $s$ and
$t$ are two maps from  $Q_1$ to $Q_0$. For any arrow $a \in Q_1$,
$s(a)$ and $t(a)$ are called its start vertex and end vertex,
respectively, and $a$ is called an arrow from $s(a)$ to $t(a)$.
For any $n\geq 0$, an $n$-path or a path of length $n$ in the
quiver $Q$ is an ordered sequence of arrows $p=a_na_{n-1}\cdots
a_1$ with $t(a_i)=s(a_{i+1})$ for all $1\leq i\leq n-1$. Note that
a 0-path is exactly a vertex and a 1-path is exactly an arrow. In
this case, we define $s(p)=s(a_1)$, the start vertex of $p$, and
$t(p)=t(a_n)$, the end vertex of $p$. For a 0-path $x$, we have
$s(x)=t(x)=x$. Let $Q_n$ be the set of $n$-paths. Let $^yQ_n^x$
denote the set of all $n$-paths from $x$ to $y$, $x, y\in Q_0$.
That is, $^yQ_n^x=\{p\in Q_n\mid s(p)=x, t(p)=y\}$.

A quiver $Q$ is {\it finite} if $Q_0$ and $Q_1$ are finite sets. A
quiver $Q$ is {\it locally finite} if $^yQ_1^x$ is a finite set
for any $x, y\in Q_0$.

Let $G$ be a group. Let ${\mathcal K}(G)$ denote the set of
conjugate classes in $G$. A formal sum $r=\sum_{C\in {\mathcal
K}(G)}r_CC$  of conjugate classes of $G$  with cardinal number
coefficients is called a {\it ramification} (or {\it ramification
data} ) of $G$, i.e.  for any $C\in{\mathcal K}(G)$, \  $r_C$ is a
cardinal number. In particular, a formal sum $r=\sum_{C\in {\mathcal
K}(G)}r_CC$  of conjugate classes of $G$ with non-negative integer
coefficients is a ramification of $G$.

 For any ramification $r$ and a $C \in {\cal K}(G)$, since $r_C$ is
 a cardinal number,
we can choice a set $I_C(r)$ such that its cardinal number is $r_C$
without lost generality.
 Let ${\mathcal K}_r(G):=\{C\in{\mathcal
K}(G)\mid r_C\not=0\}=\{C\in{\mathcal K}(G)\mid
I_C(r)\not=\emptyset\}$.  If there exists a ramification $r$ of $G$
such that the cardinal number of $^yQ_1^x$ is equal to $r_C$ for any
$x, y\in G$ with $x^{-1}y \in C\in {\mathcal K}(G)$, then $Q$ is
called a {\it Hopf quiver with respect to the ramification data
$r$}. In this case, there is a bijection from $I_C(r)$ to $^yQ_1^x$,
and hence we write  ${\ }^yQ_1^x=\{a_{y,x}^{(i)}\mid i\in I_C(r)\}$
for any $x, y\in G$ with $x^{-1}y \in C\in {\mathcal K}(G)$. Denote
by $ (Q, G, r)$ the Hopf quiver of $G$ with respect to $r$.

 The coset decomposition of $Z_{u(C)}$ in $G$ is
\begin {eqnarray} \label {e0.1}
G &=&\bigcup_{\theta\in\Theta_C}Z_{u(C)}g_{\theta},
\end {eqnarray}
where $\Theta_C$ is an index set. It is easy to check that
$|\Theta_C|=|C|$. We always assume that the representative element
of the coset $Z_{u(C)}$ is the identity $1$ of $G$. We claim that
$\theta=\eta$ if
$g_{\theta}^{-1}u(C)g_{\theta}=g_{\eta}^{-1}u(C)g_{\eta}$. In fact,
by the equation
$g_{\theta}^{-1}u(C)g_{\theta}=g_{\eta}^{-1}u(C)g_{\eta}$ one gets
$g_{\eta} g_{\theta}^{-1}u(C) = u(C)g_{\eta} g_{\theta}^{-1}$. Thus
$g_{\eta} g_{\theta}^{-1}\in Z_{u(C)}$, and hence $\theta= \eta$.
For any $x, y \in G$ with $x^{-1}y\in C\in {\cal K}(G)$, there
exists a unique $\theta\in \Theta _C$ such that \begin {eqnarray}
\label {e0.2} x^{-1}y = g_{\theta}^{-1}u(C)g_{\theta}.
\end {eqnarray}
Without specification, we will always assume that $x, y , \theta$
and $C$ satisfy the above relation (\ref{e0.2}). Note that
$\theta$ is only determined by $x^{-1}y$. For any $h\in G$ and
$\theta\in\Theta_C$, there exist unique $h'\in Z_{u(C)}$ and
$\theta'\in\Theta_C$ such that $g_{\theta}h = h'g_{\theta'}$. Let
$\zeta_{\theta}(h)=h'$. Then we have
\begin {eqnarray} \label {e0.3} g_{\theta}h&=&\zeta_{\theta}(h)g_{\theta'}.
\end {eqnarray}
If $u(C)$ lies in the center $Z(G)$ of $G$, we have
$\zeta_{\theta}=id_G$. In particular, if G is abelian, then
$\zeta_{\theta}=id_G$ since $Z_{u(C)}=G$.

Let $H$ be a Hopf algebra. A (left-left) {\it Yetter-Drinfeld
module} $V$ over $H$ (simply, {\rm YD} $H$-module) is simultaneously
a left $H$-module and a left $H$-comodule satisfying the following
compatibility condition:
\begin{eqnarray}\label{ydm}
\sum(h\cdot v)_{(-1)}\otimes(h\cdot v)_{(0)}=\sum
h_{(1)}v_{(-1)}S(h_{(3)})\otimes h_{(2)}\cdot v_{(0)},\ \ v\in V,
h\in H.
\end{eqnarray}
We denote by $^H_H{\mathcal YD}$ the category of {\rm YD}
$H$-modules; the morphisms in this category preserve both the action
and the coaction of $H$.

The structure of a Nichols algebra appeared first in the paper
\cite{Ni78},  and N. Andruskiewitsch and H. J. Schneider used it to
classify finite-dimensional pointed Hopf algebras \cite {AS98a,
AS98b, AS02, AS00}. Its definition can be found in \cite [Definition
2.1] {AS02}.

If  $\phi: A\rightarrow A'$ is  an algebra homomorphism and  $(M,
\alpha ^-)$ is a left $A'$-module, then $M$ becomes a left
$A$-module with the $A$-action given by $a \cdot x =\phi (a) \cdot x
$ for any $a\in A$, $x\in M$, called a pullback $A$-module through
$\phi$, written as  $_{\phi}M$.  Dually, if  $\phi: C\rightarrow C'$
be a coalgebra homomorphism and  $(M, \delta ^- )$ is  a left
$C$-comodule, then $M$ is a left $C'$-comodule with the
$C'$-comodule structure given by $ {\delta'}^-:=(\phi\otimes{\rm
id})\delta^-$, called  a push-out $C'$-comodule through $\phi$,
written as  $^{\phi}M$.

Let $A$ be an algebra and $M$ be an $A$-bimodule. Then the tensor
algebra $T_A(M)$ of $M$ over $A$ is a graded algebra with
$T_A(M)_{(0)}=A$, $T_A(M)_{(1)}=M$ and $T_A(M)_{(n)}=\otimes^n_AM$
for $n>1$. That is, $T_A(M)=A\oplus(\bigoplus_{n>0}\otimes^n_AM)$.
Let $D$ be another algebra. If $h$ is an algebra map from $A$ to $D$
and $f$ is an $A$-bimodule map from $M$ to $_hD_h$, then by the
universal property of $T_A(M)$ (see \cite [Proposition 1.4.1]
{Ni78}) there is a unique algebra map $T_A(h,f): T_A(M)\rightarrow
D$ such that $T_A(h,f)\iota_0=h$ and $T_A(h,f)\iota_1=f$.  One can
easily see that $T_A ( h, f ) = h + \sum _{n>0} \mu ^{n-1}T_n (f )$,
where $T_n(f)$ is the map from $\otimes _A^n M$ to $\otimes_A^nD$
given by $T_n(f)(x_1\otimes x_2 \otimes \cdots \otimes x_n) =
f(x_1)\otimes f(x_2) \otimes \cdots \otimes f(x_n)$, i.e.,
$T_n(f)=f\otimes _A f\otimes_A\cdots\otimes_A f$. Note that $\mu$
can be viewed as a map from $D\otimes _A D$ to $D$.

Dually, let $C$ be a coalgebra and let $M$ be a $C$-bicomodule. Then
the cotensor coalgebra $T_C^c(M)$ of $M$ over $C$ is a graded
coalgebra with $T_C^c(M)_{(0)}=C$, $T_C^c(M)_{(1)}=M$ and
$T_C^c(M)_n=\Box^n_CM$ for $n>1$. That is,
$T_C^c(M)=C\oplus(\bigoplus_{n>0}\Box^n_CM)$.

\section{\bf {\rm PM} quiver Hopf algebras }\label{s1}

In this section we describe {\rm PM} quiver Hopf algebras by
parameters.

We first describe the category of Hopf bimodules by categories of
modules.

Let $G$ be a group and let $(B, \delta ^-, \delta ^+)$ be a
$kG$-bicomodule. Then the $(x, y)$-isotypic component of $B$ is
$$^yB^x= \{b \in B \mid  \delta ^-(b) = y \otimes b, \delta ^+ (b)
  = b \otimes x  \},$$
where $x, y \in G$. Let $M$ be another $kG$-bicomodule and $f:
B\rightarrow M$ be a $kG$-bicomodule homomorphism. Then
$f(^yB^x)\subseteq\ ^yM^x$ for any $x, y\in G$. Denote by $^yf^x$
the restriction map $f|_{^yB^x}:\ ^yB^x\rightarrow\ ^yM^x$, $x,
y\in G$.

\begin {Theorem} \label {2} (See \cite [Proposition 3.3]{CR97} and \cite [Theorem 4.1]{CR02})
The category $^{kG}_{kG}{\mathcal M}^{kG}_{kG}$ of $kG$-Hopf bimodules
is equivalent to the Cartesian product category  $\prod_{C \in
{\mathcal K}(G)} {\mathcal M}_{kZ_{u(C)}}$ of categories ${\mathcal
M}_{kZ_{u(C)}}$  of right $kZ_{u(C)}$-modules for all $C \in
{\mathcal K}(G)$.
\end {Theorem}

For later use, we give the mutually inverse functors between the two
categories here. The functors $W$ from $^{kG}_{kG}{\mathcal
M}^{kG}_{kG}$ to $\prod_{C \in {\mathcal K}(G)}{\mathcal
M}_{kZ_{u(C)}}$ is defined by
$$W(B)=\{^{u(C)}\! B^1\}_{C\in{\mathcal K}(G)},\
W(f)=\{^{u(C)}\! f^1\}_{C\in{\mathcal K}(G)}$$ for any object $B$
and morphism $f$ in $^{kG}_{kG}{\mathcal M}^{kG}_{kG}$, where the
right $kZ_{u(C)}$-module action on ${}^{u(C)}\!  B^1$ is given by
\begin {eqnarray}\label {th1e1}
b \lhd  h = h^{-1} \cdot b \cdot h,\ \ b\in\ ^{u(C)}\! B^1,\ h\in
Z_{u(C)}.
\end {eqnarray}
The functor $V$ from $\prod_{C\in{\mathcal K}(G)}{\mathcal
M}_{kZ_{u(C)}}$ to $^{kG}_{kG}{\mathcal M}^{kG}_{kG}$ is defined as
follows:

For $M =\{ M(C)\}_{C \in {\mathcal K}(G)}\in \prod_{C \in {\mathcal
K}(G)} {\mathcal M}_{kZ_{u(C)}}$, $V(M)$ is given by
\begin {eqnarray}\label {e1.2} \begin {array}{llll}
^yV(M)^x &=& x \otimes M(C)\otimes _{kZ_{u(C)}} g_{\theta},\\
h \cdot (x \otimes m \otimes _{kZ_{u(C)}} g_{\theta})  &=&
hx \otimes m \otimes _{kZ_{u(C)}} g_{\theta},\\
 (x \otimes m \otimes _{kZ_{u(C)}} g_{\theta}) \cdot h  &=&
xh \otimes m \otimes _{kZ_{u(C)}} g_{\theta}h=xh \otimes (m\lhd
\zeta_{\theta}(h))\otimes_{kZ_{u(C)}}g_{\theta'},
\end {array}
\end {eqnarray}
where $h, x, y \in G$ with $x^{-1}y\in C$ and the relation
(\ref{e0.2}) and (\ref {e0.3}), $m\in M(C)$. For any morphism
$f=\{f_C\}_{C\in{\mathcal K}(G)}: \{M(C)\}_{C\in{\mathcal
K}(G)}\rightarrow \{N(C)\}_{C\in{\mathcal K}(G)}$, $V(f)(x\otimes
m\otimes_{kZ_{u(C)}}g_{\theta})=x\otimes
f_C(m)\otimes_{kZ_{(u(C))}}g_{\theta}$ for any $m\in M_{C}$, $x, y
\in G$ with $x^{-1}y=g^{-1}_{\theta}u(C)g_{\theta}$. That is,
$^yV(f)^x={\rm id}\otimes f_C\otimes{\rm id}$.

An $A$-module $M$ is called {\it pointed} if $M=0$ or $M$ is a
direct sum of one dimensional $A$-modules.

\begin{Definition}\label{1.1}
A $kG$-Hopf bimodule $N$ is called a {\it
 $kG$-Hopf bimodule with pointed module structure ( or a  \ {\rm PM} $kG$-Hopf bimodule in
short) } if there exists an object \ \ \  \ \ \ \ $\prod_{C \in
{\mathcal K}(G)} M(C) $ in $\prod_{C \in {\mathcal K}(G)} {\mathcal
M}_{kZ_{u(C)}}$ such that $M(C)$ is a right pointed
$kZ_{u(C)}$-module for  any $C \in {\mathcal K }(G)$ and $N\cong
V(\{ M(C)\}_{C \in {\mathcal K}(G)}):=\bigoplus _{y = g_\theta ^{-1}
u(C) g_\theta, \  x, y \in G} \ x \otimes M(C)\otimes _{kZ_{u(C)}}
g_{\theta}$ as $kG$-Hopf bimodules. Here $V(\{ M(C)\}_{C \in
{\mathcal K}(G)}):=\bigoplus _{y = g_\theta ^{-1} u(C) g_\theta, \
x, y \in G}\   x \otimes M(C)\otimes _{kZ_{u(C)}} g_{\theta}$ is
defined in the proof of Theorem \ref {2}.
\end {Definition}

Lemma \ref {1.2} -- \ref {1.6}, Theorem \ref {3} and Lemma \ref
{2.7} are well-known.

\begin{Lemma}\label{1.2}  Assume that $G$ is a finite commutative group
of exponent $m$. If $k$ contains a primitive $m$-th root of 1, then
(i) every $kG$-module is a pointed module;
 (ii)  every $kG$-Hopf
bimodule is {\rm PM}.
\end {Lemma}

\begin{Lemma}\label{1.3} Let  $H=\bigoplus_{i\geq 0}H_{(i)}$ be  a graded Hopf algebra.
Set $B:=H_{(0)}$ and $M:=H_{(1)}$. Then $M$ is a $B$-Hopf bimodule
with the $B$-actions and $B$-coactions given by
\begin {eqnarray*}\label {e2.1.3}
\alpha^-=\pi_1\mu(\iota_0\otimes\iota_1),&
&\alpha^+=\pi_1\mu(\iota_1\otimes\iota_0), \
\delta^-=(\pi_0\otimes\pi_1)\Delta\iota_1, \
 \delta^+ =(\pi_1\otimes\pi_0)\Delta\iota_1.
\end {eqnarray*}
\end {Lemma}

\begin{Lemma}\label{1.4}
{\rm (i)} \ Let $B$ be an algebra and $M$ be a $B$-bimodule. Then
the tensor algebra $T_B(M)$ of $M$ over $B$ admits a graded Hopf
algebra structure if and only if $B$ admits a Hopf algebra structure
and $M$ admits a $B$-Hopf bimodule structure.

 {\rm (ii)} Let $B$ be
a coalgebra and  $M$ be a $B$-bicomodule. Then the cotensor
coalgebra $T_B^c(M)$ of $M$ over $B$ admits a graded Hopf algebra
structure if and only if $B$ admits a Hopf algebra structure and $M$
 admits a $B$-Hopf bimodule structure.
\end {Lemma}

Let $B$ be a bialgebra (Hopf algebra) and $M$ be a $B$-Hopf
bimodule. Then $T_B^c(M)$ is a graded bialgebra (Hopf algebra) by
Lemma \ref{1.4}. Let $B[M]$ denote the subalgebra of $T_B^c(M)$
generated by $B$ and $M$. Then $B[M]$ is a bialgebra (Hopf algebra)
of type one by \cite [section 2.2, p.1533]{Ni78}.  $B[M]$ is a
graded subspace of  $T_B^c (M).$
\begin {Lemma} \label {1.5}
Let $B$ and $B'$ be two  Hopf algebras. Let $M$ and $M'$ be $B$-Hopf
bimodule and  $B'$-Hopf bimodule, respectively. Assume that $\phi:
B\rightarrow B'$ be a Hopf algebra map. If $\psi $ is simultaneously
a $B$-bimodule and $B'$-bicomodule map from $^{\phi}M^{\phi}$ to
$_{\phi}M'_{\phi}$,  then {\rm(i)} \ $T_B(\iota_0\phi,
\iota_1\psi):=\iota_0\phi+ \sum_{n>0}\mu^{n-1}T_n(\iota_1\psi)$ is a
graded Hopf algebra map from $T_B(M)$ to $T_{B'}(M')$.
 {\rm(ii)} $T_{B'}^c(\phi\pi_0,
\psi\pi_1):=\phi\pi_0+ \sum_{n>0}T_n^c (\psi\pi_1)\Delta_{n-1}$ is a
graded Hopf algebra map from $T_B^c(M)$ to $T_{B'}^c(M')$.
\end {Lemma}

\begin {Lemma} \label {1.6} Let $B$ and $B'$ be two  Hopf algebras. Let
 $M$ and $M'$ be  $B$-Hopf
bimodule and  $B'$-Hopf bimodule, respectively.   Then the following
statements are equivalent: {\rm(i)} There exists a Hopf algebra
isomorphism $\phi: B\rightarrow B'$ such that $M\cong\ _{\phi}
^{\phi ^{-1}}M'{}_{\phi} ^{\phi ^{-1}}$ as $B$-Hopf bimodules.
 {\rm(ii)} $T_{B}(M)$ and $T_{B'}(M')$
are isomorphic as graded Hopf algebras.
 {\rm(iii)} $T_{B}^c(M)$ and
$T_{B'}^c(M')$ are isomorphic as graded Hopf algebras.
 {\rm(iv)} $B[M]$ and
$B'[M']$ are isomorphic as graded Hopf algebras.
\end{Lemma}

Let $Q=(G, Q_1, s, t)$ be a quiver of a group $G$. Then $kQ_1$
becomes  a $kG$-bicomodule  under  the natural comodule structures:
\begin{eqnarray}\label{arcom}
\delta^-(a)=t(a)\otimes a,\ \ \delta^+(a)=a\otimes s(a),\ \ a\in
Q_1,
\end{eqnarray} called an {\it arrow comodule}, written as $kQ_1 ^c$. In this case,
the path coalgebra $kQ^c$ is exactly isomorphic to the cotensor
coalgebra $T^c_{kG}(kQ_1^c)$  over $kG$ in a natural way (see
\cite{CM97} and \cite{CR02}). We will regard $kQ^c=T^c_{kG}(kQ_1^c)$
in the following. Moreover, when $G$ is finite,  $kQ_1$ becomes a
$(kG)^*$-bimodule with the module structures defined by
\begin{eqnarray}\label{armod}
\mbox{\hspace{1cm}}p\cdot a:=\langle p, t(a)\rangle a,\ \ a\cdot
p:=\langle p, s(a)\rangle a,\ \ p\in(kG)^*, a\in Q_1,
\end {eqnarray}
written as $kQ_1^a$, called an {\it arrow module}. Therefore, we
have a tensor algebra $T_{(kG)^*}(kQ_1^a)$. Note that the tensor
algebra $T_{(kG)^*}(kQ_1^a)$ of $kQ_1^a$ over $(kG)^*$ is exactly
isomorphic to the path algebra $kQ^a$. We will regard
$kQ^a=T_{(kG)^*}(kQ_1^a)$ in the following.

 Assume that $Q$ is a finite  quiver on finite
group $G$. Let $\xi _{kQ_1^a}$ denote the linear map from $kQ_1^a$
to $(kQ_1^c)^*$ by sending $a$ to $a^*$ for  any $a\in Q_1$ and $\xi
_{kQ_1^c}$ denote the linear map from $kQ_1^c$ to $(kQ_1^a)^*$ by
sending $a$ to $a^*$ for  any $a\in Q_1$. Here $\{ a^* \mid a^*  \in
(kQ_1)^*\}$ is the dual basis of $\{ a \mid a \in Q_1 \}$.

\begin {Lemma} \label {1.7} Assume that $Q$ is a finite Hopf quiver on finite
group $G$. Then

(i) If $(M, \alpha ^-, \alpha ^+, \delta^-, \delta ^+)$ is a finite
dimensional $B$-Hopf bimodule and $B$ is a finite dimensional Hopf
algebra, then $(M^*, \delta^{-*}, \delta ^{+*}, \alpha ^{-*}, \alpha
^{+*})$ is a $B^*$-Hopf bimodule.

(ii) If $(kQ_1^c, \alpha ^-, \alpha ^+, \delta^-, \delta ^+)$ is a
$kG$- Hopf bimodule, then there exist  unique left $(kG)^*$-comodule
operation $ \delta _{kQ_1^a}^-$ and right $(kG)^*$-comodule
operation
 $\delta
_{kQ_1^a}^+$ such that $(kQ_1^a, \alpha _{kQ_1^a} ^-, \alpha
_{kQ_1^a}^+,$ \ $ \delta_{kQ_1^a}^-, $ \ $ \delta_{kQ_1^a} ^+)$
becomes a $(kG)^*$-Hopf bimodule and $\xi_{kQ_1^a}$ becomes a
$(kG)^*$-Hopf bimodule isomorphism from $(kQ_1^a, \alpha _{kQ_1^a}
^-, \alpha _{kQ_1^a}^+, \delta_{kQ_1^a}^-, \delta_{kQ_1^a} ^+)$ to
$((kQ_1^c)^*, \delta^-{}^*, \delta ^+{}^*, \alpha ^-{}^*, \alpha
^+{}^* )$.

(iii) If $(kQ_1^a, \alpha ^-, \alpha ^+, \delta^-, \delta ^+)$ is a
$(kG)^*$- Hopf bimodule, then there exist  unique left $kG$-module
operation $ \alpha _{kQ_1^c}^-$ and right $kG$-module
 operation $\alpha
_{kQ_1^c}^+$ such that $(kQ_1^c, \alpha _{kQ_1^c} ^-, \alpha
_{kQ_1^c}^+, \delta_{kQ_1^c}^-,$ $ \delta_{kQ_1^c} ^+)$ become a
$kG$-Hopf bimodule and $\xi_{kQ_1^c}$ becomes a $kG$-Hopf bimodule
isomorphism from $(kQ_1^c, \alpha _{kQ_1^c} ^-, \alpha _{kQ_1^c}^+,
\delta_{kQ_1^c}^-, \delta_{kQ_1^c} ^+)$ to $((kQ_1^a)^*,
\delta^-{}^*, \delta ^+{}^*, \alpha ^-{}^*, \alpha ^+{}^* )$.

(iv)
 $\xi_{kQ_1^a}$ is a $(kG)^*$-Hopf
 bimodule isomorphism from $(kQ_1^a, \alpha _{kQ_1^a} ^-, \alpha
 _{kQ_1^a}^+, \delta_{kQ_1^a}^-, \delta_{kQ_1^a} ^+)$ to
 $((kQ_1^c)^*, \delta_{kQ_1^c}^-{}^*, \delta_{kQ_1^c} ^+{}^*, \alpha_{kQ_1^c} ^-{}^*, \alpha _{kQ_1^c}
 ^+{}^* )$ if and only if
 $\xi_{kQ_1^c}$ becomes a $kG$-Hopf bimodule
 isomorphism from $(kQ_1^c, \alpha _{kQ_1^c} ^-, \alpha _{kQ_1^c}^+,
 \delta_{kQ_1^c}^-, \delta_{kQ_1^c} ^+)$ to $((kQ_1^a)^*,$ $
 \delta_{kQ_1^a}^-{}^*, $ $ \delta_{kQ_1^a} ^+{}^*, $ $\alpha _{kQ_1^a}^-{}^*, \alpha _{kQ_1^a}^+{}^* )$.
\end {Lemma}

{\bf Proof.} It is easy to check (i)--(iii). Now we show (iv). Let
$A:= kQ^a_1,$ $ B:= kQ_1^c.$  Let  $\sigma_ B$ denote the canonical
linear isomorphism from  $B$ to $ B^{**}$ by sending $b$ to $b^{**}$
for any $b \in B$, where  $<b^{**}, f> = <f, b>$ for any $f \in
B^*.$ If $A \stackrel {\xi _A } {\cong } B^* $ as $(kG)^*$-Hopf
bimodules, then $B \stackrel {\sigma _B} {\cong } B^{**}  \stackrel
{(\xi _A )^*} {\cong }A^*$ as $kG$-Hopf bimodules. It is easy to
check $\xi _B =(\xi _A)^* \sigma _B $. Therefore $\xi _B$ is a
$kG$-Hopf bimodule isomorphism. Conversely, if $B \stackrel {\xi _B
} {\cong } A^* $ as $kG$-Hopf bimodules, we can similarly show that
$A \cong  B^* $ as $(kG)^*$-Hopf bimodules. \ $\Box$

\begin {Theorem} \label {3} (see \cite [Theorem 3.3]{CR02} and \cite [Theorem 3.1] {CR97})
Let $Q$ be a quiver over  group $G$. Then the following two
statements are equivalent:

(i) $Q$ is a Hopf quiver.

(ii) Arrow comodule $kQ_1^c$ admits a $kG$-Hopf bimodule
structure.\\
Furthermore, if $Q$ is finite, then the above are equivalent to the
following:

(iii) Arrow module $kQ_1^a$ admits a $(kG)^*$-Hopf bimodule
structure.
\end {Theorem}

Assume that $Q$ is a Hopf quiver. It follows from Theorem \ref {3}
that there exist a left $kG$-module structure $\alpha ^-$ and a
right $kG$-module structure $\alpha ^+$ on arrow comodule $(kQ_1^c,
\delta^-, \delta ^+)$ such that  $(kQ_1^c, \alpha ^-, \alpha ^+,
\delta^-, \delta ^+)$ becomes a $kG$-Hopf bimodule, called a
$kG$-Hopf bimodule with arrow comodule, written  $(kQ_1^c, \alpha
^-, \alpha ^+ )$ in short. We obtain three graded Hopf algebras
$T_{kG} (kQ_1^c)$, $T_{kG}^c(kQ_1^c)$ and $kG[kQ_1^c]$, called
semi-path Hopf algebra, co-path Hopf algebra and one-type-co-path
Hopf algebra, written $kQ^{s}(\alpha ^-, \alpha ^+)$, $kQ^{c}(
\alpha ^-, \alpha ^+)$ and $kG[kQ_1^c, \alpha ^-, \alpha ^+]$,
respectively. Dually, when $Q$ is finite, it follows from Theorem
\ref {3} that there exist a left $(kG)^*$-comodule structure $\delta
^-$ and a right $(kG)^*$-comodule structure $\delta ^+$ on arrow
module $(kQ_1^a, \alpha^-, \alpha ^+)$ such that  $(kQ_1^a, \alpha
^-, \alpha ^+, \delta^-, \delta ^+)$ becomes a $(kG)^*$-Hopf
bimodule, called a $(kG)^*$-Hopf bimodule with arrow module, written
$(kQ_1^a, \delta^-, \delta ^+)$ in short. We obtain three  graded
Hopf algebras $T_{(kG)^*} (kQ_1^a)$, $T_{(kG)^*}^c(kQ_1^a)$ and
$(kG)^* [kQ_1^a]$, called path Hopf algebra, semi-co-path Hopf
algebra and one-type-path Hopf algebra, written $kQ^a( \delta ^-,
\delta ^+)$, $kQ^{sc}( \delta ^-, \delta ^+)$ and $(kG)^*[kQ_1^a
,\alpha _1, \alpha ^+]$, respectively.  We call the six graded Hopf
algebras the quiver Hopf algebras (over $G$). We usually omit the
(co)module operations when we write these quiver Hopf algebras.

If $\xi _{kQ_1^a}$ or $\xi _{kQ_1^c}$ is a Hopf bimodule
isomorphism, then, by Lemma \ref {1.6} and Lemma \ref {1.7},
$T_{(kG)^*} (\iota _0,\iota_1 \xi _{kQ_1^a})$ and $T_{kG} ^c (\pi
_0, \xi _{kQ_1^c} \pi_1)$ are graded Hopf algebra isomorphisms from
$T_{(kG)^*} (kQ_1^a)$ to $T_{(kG)^*} ((kQ_1^c)^*)$ and from $T_{kG}
^c(kQ_1^c)$ to $T_{kG}^c((kQ_1^a)^*)$, respectively; $T_{(kG)^*} ^c
(\pi _0,$ $ \xi _{kQ_1^a} \pi_1)$ and $T_{kG}  (\iota _0, \iota _1
\xi _{kQ_1^c})$ are graded Hopf algebra isomorphisms from
$T_{(kG)^*} ^c(kQ_1^a)$ to $T_{(kG)^*} ^c((kQ_1^c)^*)$ and from
$T_{kG} (kQ_1^c)$ to $T_{kG}((kQ_1^a)^*)$, respectively. In this
case, $(kQ_1^a, kQ_1^c)$, $(kQ^a, kQ^c)$ and $(kQ^{s}, kQ^{sc})$ are
said to be arrow dual pairings.

If $(kQ_1^c,\alpha^-,\alpha^+,\delta^-,\delta^+)$ is a {\rm PM}
$kG$-Hopf bimodule, and $(kQ_1^c, kQ_1^a)$ is an arrow pairing, then
$(kQ_1^a, \alpha ^{-*}, \alpha ^{+*}) $ is called a {\rm PM}
$(kG)^*$-Hopf bimodule and  six quiver Hopf algebras induced by
$kQ_1^c$ and $kQ_1^a$ are called {\rm PM} quiver Hopf algebras.

Now we are going to describe the structure of all {\rm PM} $kG$-Hopf
bimodules and the corresponding graded Hopf algebras.

\begin{Definition}\label {1.8}
$(G, r, \overrightarrow \chi, u)$ is called a ramification system
with characters   (or {\rm RSC } in short ), if $r$ is a
ramification of $G$, $u$ is a map from ${\mathcal K}(G)$ to $G$ with
$u(C)\in C$ for any $C\in {\mathcal K}(G)$, and $\overrightarrow
\chi=\{\chi_C^{(i)} \}_ { i\in I_C(r), C\in{\mathcal K}_r(G)} \ \in
\prod _ { C\in{\mathcal K}_r(G)} (\widehat{Z_{u(C)}}) ^{r_C}$ with
$\chi_C^{(i)} \in \widehat{Z_{u(C)}} $ for any
 $ i \in I_C(r), C\in {\mathcal
K}_r(G)$.

${\rm RSC} (G, r, \overrightarrow \chi, u)$ and ${\rm RSC} (G', r',
\overrightarrow {\chi'}, u')$ are said to be {\it isomorphic} if the
following conditions are satisfied:

$\bullet$ There exists a group isomorphism $\phi: G\rightarrow G'$.

$\bullet$ For any $C\in{\mathcal K}(G)$, there exists an element
$h_C\in G$ such that $\phi(h_C^{-1}u(C)h_C)=u'(\phi(C))$.

$\bullet$ For any $C\in{\mathcal K}_r(G)$, there exists a bijective
map $\phi_C : I_C(r) \rightarrow I_{\phi(C)}(r')$ such that $\chi
'{}_{\phi(C)}^{(\phi_C(i))}(\phi(h_C^{-1}hh_C))=\chi_C^{(i)}(h)$ for
all $h\in Z_{u(C)}$ and $i\in I_C(r)$.
\end {Definition}

{\bf Remark.} Assume that $G=G'$, $r=r'$ and $u(C) = u'(C)$ for any
$C \in {\cal K}_r(G)$. If there is a permutation $\phi_C$ on
$I_C(r)$ for any $C\in{\mathcal K}_r(G)$ such that
${\chi'}_C^{(\phi_C(i))}=\chi_C^{(i)}$ for all $i\in I_C(r)$, then
obviously
 ${\rm RSC} (G, r, \overrightarrow \chi, u)\cong {\rm RSC} (G, r, \overrightarrow{\chi'},
 u)$.

\begin{Proposition}\label {1.10} If $N$ is a {\rm PM} $kG$-Hopf bimodule, then
there exist a Hopf quiver $(Q, G,r )$, an ${\rm RSC} (G, r,
\overrightarrow \chi, u)$ and a $kG$-Hopf bimodule $(kQ_1^c, \alpha
^-, \alpha ^+)$ with
$$ \alpha ^- (h \otimes a^{(i)}_{y,x}) :=   h\cdot a^{(i)}_{y,x}=a^{(i)}_{hy,hx},\ \
\alpha ^+  ( a^{(i)}_{y,x}\otimes h) := a^{(i)}_{y,x}\cdot
h=\chi^{(i)}_C(\zeta_{\theta}(h))a^{(i)}_{yh,xh}$$ where $x, y, h\in
G$ with $x^{-1}y=g^{-1}_{\theta}u(C)g_{\theta}$, $\zeta_{\theta}$ is
given by {\rm(\ref{e0.3})}, $C\in{\mathcal K}_r(G)$ and $i\in
I_C(r)$,  such that $N \cong (kQ_1^c, \alpha ^-, \alpha ^+)$ as
$kG$-Hopf bimodules.

\end {Proposition}
{\bf Proof.} Since $N$  is a {\rm PM} $kG$-Hopf bimodule, there
exists an object \ \ \  \ \ \ \ $\prod_{C \in {\mathcal K}(G)} M(C)
$ in $\prod_{C \in {\mathcal K}(G)} {\mathcal M}_{kZ_{u(C)}}$ such
that $M(C)$ is a pointed $kZ_{u(C)}$-module for  any $C \in
{\mathcal K }(G)$ and $N \cong V(\{ M(C)\}_{C \in {\mathcal
K}(G)})=\bigoplus _{y = g_\theta ^{-1} u(C) g_\theta, \  x, y \in
G}\   x \otimes M(C)\otimes _{kZ_{u(C)}} g_{\theta}$ as $kG$-Hopf
bimodules. Let $r = \sum _{C \in {\mathcal K} (G)}r_C C$ with $ r_C
={\rm dim}M(C)$ for any $C \in {\mathcal K} (G)$. Notice that ${\rm
dim}M(C)$ denotes the  cardinal number  of a basis  of a basis of
$M(C)$ when $M(C)$ is infinite dimensional.  Since $(M(C), \alpha
_C)$ is a pointed $kZ_{u(C)}$-module, there exist a $k$-basis
$\{x_C^{(i)}\mid i\in I_C(r)\}$ in $M(C)$ and a family of characters
$\{\chi_C^{(i)}\in\widehat{Z_{u(C)}}\mid i\in I_C(r)\}$ such that $
\alpha _C (x_C^{(i)} \otimes  h)= x_C^{(i)}\lhd
h=\chi_C^{(i)}(h)x_C^{(i)}$ for any $i\in I_C(r)$ and $h\in
Z_{u(C)}$.

We have to show that $(kQ_1^c, \alpha ^-, \alpha ^+)$ is isomorphic
to $\bigoplus _{y = g_\theta ^{-1} u(C) g_\theta, \  x, y \in G}\
x \otimes M(C)\otimes _{kZ_{u(C)}} g_{\theta}$ as $kG$-Hopf
bimodules. Observe that there is a canonical $kG$-bicomodule
isomorphism $\varphi: kQ_1\rightarrow \bigoplus _{y = g_\theta ^{-1}
u(C) g_\theta, \  x, y \in G}\   x \otimes M(C)\otimes _{kZ_{u(C)}}
g_\theta$ given by
\begin{eqnarray}\label{indbya}
\varphi(a^{(i)}_{y,x})=x\otimes
x_C^{(i)}\otimes_{kZ_{u(C)}}g_{\theta}
\end{eqnarray}
where $x, y\in G$ with $x^{-1}y=g^{-1}_{\theta}u(C)g_{\theta}$,
$C\in{\mathcal K}_r(G)$ and $i\in I_C(r)$. Now we have

\begin{eqnarray*}
\varphi (\alpha ^- (h \otimes a^{(i)}_{y,x})) &=& \varphi ( a^{(i)}_{hy,hx})=    hx \otimes x^{(i)}_C \otimes_{kZ_{u(C)}}g_{\theta}  \ \ \\
&=&    h \cdot (x \otimes x^{(i)}_C \otimes_{kZ_{u(C)}}g_{\theta})  \ \ \ (\hbox {see }( \ref {e1.2} )) \\
&=& h \cdot  \varphi ( a^{(i)}_{y,x})
\end{eqnarray*}
and
\begin{eqnarray*}
\varphi (\alpha ^+ (a^{(i)}_{y,x} \otimes h))
&=&\chi^{(i)}_C(\zeta_{\theta}(h))\varphi ( a^{(i)}_{yh,xh})
= \chi^{(i)}_C(\zeta_{\theta}(h))   xh \otimes x^{(i)}_C \otimes_{kZ_{u(C)}}g_{\theta'} \\
&{}&  (\hbox {since } h g_{\theta'}^{-1}=g_\theta ^{-1} \zeta
_\theta (h)
 \hbox { and }
yh = x h g_{\theta '}^{-1}u(C)g_{\theta '})\\
&=&   (x \otimes x^{(i)}_C \otimes_{kZ_{u(C)}}g_{\theta}) \cdot h\ \ \ (\hbox { see } (\ref {e1.2} )) \\
&=&  \varphi ( a^{(i)}_{y,x}) \cdot h,
\end{eqnarray*}
 where $x, y, h\in
G$ with $x^{-1}y=g^{-1}_{\theta}u(C)g_{\theta}$, $\zeta_{\theta}$ is
given by {\rm(\ref{e0.3})}, $C\in{\mathcal K}_r(G)$ and $i\in
I_C(r)$. Consequently, $\varphi$ is a $kG$-Hopf bimodule
isomorphism. $\Box$

Let $(kQ_1^c, G, r, \overrightarrow \chi, u)$ denote the $kG$-Hopf
bimodule $(kQ_1^c, \alpha ^-, \alpha ^+)$ given in Lemma \ref
{1.10}. Furthermore, if $(kQ_1^c, kQ_1^a)$ is an arrow dual pairing,
then we denote the $(kG)^*$-Hopf bimodule $kQ_1^a $ by $(kQ_1^a, G,
r, \overrightarrow \chi, u)$. We obtain six quiver Hopf algebras
$kQ^c (G, r, \overrightarrow \chi, u)$, $kQ^s (G, r, \overrightarrow
\chi, u)$, $kG[ kQ_1^c, G, r, \overrightarrow \chi, u]$, $kQ^a (G,
r, \overrightarrow \chi, u)$, $kQ^{sc} (G, r, \overrightarrow \chi,
u),$ $(kG)^*[ kQ_1^a, G, r, \overrightarrow \chi,$ $ u]$, called the
quiver Hopf algebras determined by ${\rm RSC} (G,$ $ r,$ $
\overrightarrow \chi,u)$.

From Proposition \ref{1.10}, it seems that the right $kG$-action on
$(kQ_1^c,G, r, \overrightarrow \chi, u)$ depends on the choice of
the set $\{g_{\theta}\mid \theta\in\Theta_C\}$ of coset
representatives of $Z_{u(C)}$ in $G$ (see, Eq.(\ref{e0.1})). The
following lemma shows that $(kQ_1^c,G, r, \overrightarrow \chi, u)$
is, in fact, independent of the choice of the coset representative
set $\{g_{\theta}\mid \theta\in\Theta_C\}$, up to $kG$- Hopf
bimodule isomorphisms. For a while, we write $(kQ_1^c,G, r,
\overrightarrow \chi, u)=(kQ_1^c,G, r, \overrightarrow \chi, u,
\{g_{\theta}\})$ given before. Now let $\{h_{\theta}\in G\mid
\theta\in\Theta_C\}$ be another coset representative set of
$Z_{u(C)}$ in $G$ for any $C\in{\mathcal K}(G)$. That is,
\begin{eqnarray}\label{ncosetde}
G=\bigcup_{\theta\in\Theta_C}Z_{u(C)}h_{\theta}.
\end{eqnarray}
\begin{Lemma}\label{1.13}
With the above notations, $(kQ_1^c,G, r, \overrightarrow \chi, u,
\{g_{\theta}\})$ and $(kQ_1^c,G, r, \overrightarrow \chi, u, $
$\{h_{\theta}\})$ are isomorphic $kG$-Hopf bimodules.
\end{Lemma}

{\bf Proof.} We may assume \ \
$Z_{u(C)}h_{\theta}=Z_{u(C)}g_{\theta}$ \ \ for any \ $C\in{\mathcal
K}(G)$ and \ \  $\theta\in\Theta_C$. Then \ \ \
$g_{\theta}h^{-1}_{\theta}\in Z_{u(C)}$. Now let \ $x, y, h\in G$
with $x^{-1}y=g^{-1}_{\theta}u(C)g_{\theta}$. Then $x^{-1}y =$\ $
h^{-1}_{\theta}(g_{\theta}h^{-1}_{\theta})^{-1}u(C)(g_{\theta}h^{-1}_{\theta})h_{\theta}
$ \ $=h^{-1}_{\theta}u(C)h_{\theta}$ and
$h_{\theta}h=(h_{\theta}g^{-1}_{\theta})g_{\theta}h =
(h_{\theta}g^{-1}_{\theta})\zeta_{\theta}(h)g_{\theta'}=
(h_{\theta}g^{-1}_{\theta})\zeta_{\theta}(h)(g_{\theta'}h^{-1}_{\theta'})h_{\theta'}$,
where $g_{\theta}h=\zeta_{\theta}(h)g_{\theta'}$. Hence from
Proposition \ref{1.10} we know that the right $kG$-action on
$(kQ_1^c,G, r, \overrightarrow \chi, u, \{h_{\theta}\})$ is given by
$$\begin{array}{rcl}
a^{(i)}_{y,x}\cdot h& =&
\chi^{(i)}_C((h_{\theta}g^{-1}_{\theta})\zeta_{\theta}(h)(g_{\theta'}h^{-1}_{\theta'}))
a^{(i)}_{yh,xh}\\
&=&
\chi^{(i)}_C((g_{\theta}h^{-1}_{\theta})^{-1})\chi^{(i)}_C(\zeta_{\theta}(h))
\chi^{(i)}_C(g_{\theta'}h_{\theta'}^{-1}) a^{(i)}_{yh,xh},\ i\in
I_C(r).
\end{array}$$
However, we also have
$$(xh)^{-1}(yh)=h^{-1}(x^{-1}y)h=h^{-1}g^{-1}_{\theta}u(C)g_{\theta}h=g^{-1}_{\theta'}u(C)g_{\theta'}.$$
It follows that the $k$-linear isomorphism $f: kQ_1\rightarrow
kQ_1$ given by
$$f(a^{(i)}_{y,x})=\chi^{(i)}_C(g_{\theta}h^{-1}_{\theta})a^{(i)}_{y,x}$$
for any $x,y\in G$ with $x^{-1}y=g^{-1}_{\theta}u(C)g_{\theta}$,
$C\in{\mathcal K}_r(G)$ and $i\in I_C(r)$, is a $kG$-Hopf bimodule
isomorphism from $(kQ_1^c,G, r, \overrightarrow \chi, u,
\{g_{\theta}\})$ to $(kQ_1^c,G, r, \overrightarrow \chi, u,
\{h_{\theta}\})$. \  \ $\Box$

Now we state one of our main results, which classifies the {\rm PM}
(co-)path Hopf algebras, {\rm PM} semi-(co-)path Hopf algebras and
{\rm PM} one-type-path Hopf algebras.

\begin {Theorem} \label {4}
  Let $(G, r, \overrightarrow{\chi}, u)$ and
$(G', r', \overrightarrow{\chi'}, u')$ are two {\rm RSC}'s. Then the
following statements are equivalent:

{\rm(i)} ${\rm RSC} (G, r , \overrightarrow{\chi}, u)$ $\cong $
${\rm RSC} (G', r', \overrightarrow{\chi'}, u')$.

{\rm(ii)} There exists a Hopf algebra isomorphism $\phi:
kG\rightarrow kG'$ such that $(kQ_1^c, G, r , \overrightarrow{\chi},
u)  \cong\ _{\phi} ^{\phi ^{-1}}
 ( ( kQ_1' {}^c, G', r', \overrightarrow{\chi'}, u') ){}_{\phi} ^{\phi ^{-1}}$ as $kG$-Hopf bimodules.

{\rm(iii)} $kQ^c(G, r, \overrightarrow \chi, u)\cong k{Q'}^c(G, r,
\overrightarrow \chi, u)$.
 {\rm(iv)}  $kQ^s(G, r, \overrightarrow \chi, u) \cong
k{Q'}^s(G', r', \overrightarrow \chi', u')$.

 {\rm(v)} \
$kG[kQ_1^c,G, r, \overrightarrow \chi, u] $ $\cong kG'[kQ_1'{}^c,G', r', \overrightarrow \chi', u']$.\\
Furthermore, if $Q$ is finite, then the above are equivalent to the
following:

{\rm(vi)} $kQ^a(G, r, \overrightarrow \chi, u) $ $ \cong k{Q'}^a(G',
r', \overrightarrow \chi', u')$. {\rm(vii)} $kQ^{sc}(G, r,
\overrightarrow \chi, u) \cong k{Q'}^{sc}(G', r', \overrightarrow
\chi', u')$. {\rm(viii)} \ $(kG)^*[kQ_1^a, G, r, \overrightarrow
\chi, u]\cong (kG')^*[kQ_1'{}^a, G', r', \overrightarrow \chi',
u']$. Notice that the isomorphisms above are ones of graded Hopf
algebras but (i) (ii).
\end {Theorem}

{\bf Proof.}  By Lemma \ref{1.6} and Lemma \ref {1.7}, we only have
to prove (i) $\Leftrightarrow$ (ii).

(i) $\Rightarrow$ (ii). Assume that ${\rm RSC} (G, r ,
\overrightarrow{\chi}, u)$ $\cong $ ${\rm RSC} (G', r',
\overrightarrow{\chi'}, u')$. Then  there exist a group isomorphism
$\phi: G\rightarrow G'$,  an element $h_C\in G$ such that
$\phi(h^{-1}_Cu(C)h_C)=u'(\phi(C))$ for any $C\in{\mathcal K}(G)$
and a bijective map $\phi_C: I_C(r)\rightarrow I_{\phi(C)}(r')$ such
that
${\chi'}_{\phi(C)}^{(\phi_C(i))}(\phi(h^{-1}_Chh_C))=\chi_C^{(i)}(h)$
for any $h\in Z_{u(C)}$, $C\in{\mathcal K}_r(G)$ and $i\in I_C(r)$.
Then $\phi(h^{-1}_CZ_{u(C)}h_C)=Z_{u'(\phi(C))}$ and $\phi:
kG\rightarrow kG'$ is a Hopf algebra isomorphism. Now let
$G=\bigcup_{\theta\in\Theta_C}Z_{u(C)}g_{\theta}$ be given as in
(\ref{e0.2}) for any $C\in{\mathcal K}(G)$, and assume that the
$kG$-Hopf bimodule $(kQ_1^c, G, r, \overrightarrow \chi, u)$ is
defined by using these coset decompositions. Then
\begin{eqnarray}\label{coset}
G'=\bigcup_{\theta\in\Theta_C}Z_{u'(\phi(C))}(\phi(h^{-1}_Cg_{\theta}h_C))
\end{eqnarray}
is a coset decomposition of $Z_{u'(\phi(C))}$ in $G'$ for any
$\phi(C)\in{\mathcal K}(G')$. By Lemma \ref{1.13}, we may assume
that the structure of the $kG'$-Hopf bimodule $(kQ_1'{}^c, G', r',
\overrightarrow \chi', u')$ is obtained by using these coset
decompositions (\ref{coset}). Define a $k$-linear isomorphism $\psi:
(kQ_1, G, r, \overrightarrow \chi, u) \rightarrow (kQ_1', G', r',
\overrightarrow \chi', u')$ by
$$\psi(a_{y,x}^{(i)})=\chi_C^{(i)}(\zeta_{\theta}(h^{-1}_C))a_{\phi(y),\phi(x)}^{(\phi_C(i))}$$
for any $x, y\in G$ with $x^{-1}y=g^{-1}_{\theta}u(C)g_{\theta}$,
and $i\in I_C(r)$, where $C\in{\mathcal K}_r(G)$ and
$g_{\theta}h^{-1}_C=\zeta_{\theta}(h^{-1}_C)g_{\eta}$ with
$\zeta_{\theta}(h^{-1}_C)\in Z_{u(C)}$ and $\theta, \eta\in
\Theta_C$. It is easy to see that $\psi$ is a $kG$-bicomodule
homomorphism from $(kQ_1^c, G, r, \overrightarrow \chi, u)$ to
$^{\phi^{-1}}_{\phi}(kQ_1'{}^c, G', r', \overrightarrow \chi',
u')^{\phi^{-1}}_{\phi}$. Since $(hx)^{-1}(hy)=x^{-1}y$ for any $x,
y, h\in G$, it follows from Proposition \ref{1.10} that $\psi$ is
also a left $kG$-module homomorphism from $(kQ_1^c, G, r,
\overrightarrow \chi, u)$ to $^{\phi^{-1}}_{\phi}(kQ_1'{}^c, G, r,
\overrightarrow \chi, u)^{\phi^{-1}}_{\phi}$.

Now let $x, y, h\in G$ with $x^{-1}y=g^{-1}_{\theta}u(C)g_{\theta}$,
$C\in{\mathcal K}_r(G)$ and $\theta\in\Theta_C$. Assume that
$g_{\theta}h^{-1}_C=\zeta_{\theta}(h^{-1}_C)g_{\eta}$,
$g_{\theta}h=\zeta_{\theta}(h)g_{\theta'}$,
$g_{\eta}(h_Chh^{-1}_C)=\zeta_{\eta}(h_Chh^{-1}_C)g_{\eta'}$ and
$g_{\theta'}h^{-1}_C=\zeta_{\theta'}(h^{-1}_C)g_{\theta''}$ with
$\zeta_{\theta}(h^{-1}_C), \zeta_{\theta}(h),
\zeta_{\eta}(h_Chh^{-1}_C), \zeta_{\theta'}(h^{-1}_C)\in Z_{u(C)}$
and $\eta, \theta', \eta', \theta''\in \Theta_C$. Then we have
$$\begin{array}{rcl}
g_{\theta}hh^{-1}_C&=&\zeta_{\theta}(h)g_{\theta'}h^{-1}_C\ = \zeta_{\theta}(h)\zeta_{\theta'}(h^{-1}_C)g_{\theta''}\\
\end{array}$$
and
$$\begin{array}{rcl}
g_{\theta}hh^{-1}_C&=&(g_{\theta}h^{-1}_C)(h_Chh^{-1}_C)=\zeta_{\theta}(h^{-1}_C)g_{\eta}(h_Chh^{-1}_C)\
=\zeta_{\theta}(h^{-1}_C)\zeta_{\eta}(h_Chh^{-1}_C)g_{\eta'}.
\end{array}$$
It follows that $\theta''=\eta'$ and
\begin{eqnarray}\label{coeff}
\zeta_{\theta}(h)\zeta_{\theta'}(h^{-1}_C)=
\zeta_{\theta}(h^{-1}_C)\zeta_{\eta}(h_Chh^{-1}_C).
\end{eqnarray}
By Proposition \ref{1.10}, $a^{(i)}_{y,x}\cdot
h=\chi^{(i)}_C(\zeta_{\theta}(h))a^{(i)}_{yh,xh}$ for any $i\in
I_C(r)$. Moreover, we have
$(xh)^{-1}(yh)=h^{-1}g^{-1}_{\theta}u(C)g_{\theta}h=
g^{-1}_{\theta'}u(C)g_{\theta'}$. This implies that
$$\psi(a^{(i)}_{y,x}\cdot
h)=\chi^{(i)}_C(\zeta_{\theta}(h))\psi(a^{(i)}_{yh,xh})
=\chi^{(i)}_C(\zeta_{\theta}(h))\chi^{(i)}_C(\zeta_{\theta'}(h^{-1}_C))
a^{(\phi_C(i))}_{\phi(yh),\phi(xh)}.$$ On the other hand, we have
$g_{\theta}=g_{\theta}h^{-1}_Ch_C=\zeta_{\theta}(h^{-1}_C)g_{\eta}h_C$,
and hence
$$\begin{array}{rcl}
\phi(x)^{-1}\phi(y)&=&\phi(x^{-1}y)=\phi(g^{-1}_{\theta}u(C)g_{\theta})
=\phi(h^{-1}_Cg^{-1}_{\eta}u(C)g_{\eta}h_C)\\
&=& \phi(h^{-1}_Cg^{-1}_{\eta}h_C)\phi(h^{-1}_Cu(C)h_C)\phi(h^{-1}_Cg_{\eta}h_C)\\
&=&\phi(h^{-1}_Cg_{\eta}h_C)^{-1}u'(\phi(C))\phi(h^{-1}_Cg_{\eta}h_C).
\end{array}$$
We also have
$$\begin{array}{rcl}
\phi(h^{-1}_Cg_{\eta}h_C)\phi(h)&=&\phi(h^{-1}_Cg_{\eta}h_Chh^{-1}_Ch_C)\\
&=&\phi(h^{-1}_C\zeta_{\eta}(h_Chh^{-1}_C)g_{\eta'}h_C)\\
&=&\phi(h^{-1}_C\zeta_{\eta}(h_Chh^{-1}_C)h_C)\phi(h^{-1}_Cg_{\eta'}h_C).
\end{array}$$
Thus by Proposition \ref{1.10} one gets
$$\begin{array}{rcl}
a^{(\phi_C(i))}_{\phi(y),\phi(x)}\cdot\phi(h)&=&
{\chi'}^{(\phi_C(i))}_{\phi(C)}(\phi(h^{-1}_C\zeta_{\eta}(h_Chh^{-1}_C)h_C))
a^{(\phi_C(i))}_{\phi(yh),\phi(xh)}\\
 &=&
\chi^{(i)}_C(\zeta_{\eta}(h_Chh^{-1}_C))
a^{(\phi_C(i))}_{\phi(yh),\phi(xh)}.\\
\end{array}$$
Now it follows from Eq.(\ref{coeff}) that
$$\begin{array}{rcl}
\psi(a^{(i)}_{y,x})\cdot\phi(h)&=&
\chi_C^{(i)}(\zeta_{\theta}(h^{-1}_C))a_{\phi(y),\phi(x)}^{(\phi_C(i))}
\cdot\phi(h)\\
&=&\chi_C^{(i)}(\zeta_{\theta}(h^{-1}_C))
\chi^{(i)}_C(\zeta_{\eta}(h_Chh^{-1}_C))
a^{(\phi_C(i))}_{\phi(yh),\phi(xh)}\\
&=&\psi(a^{(i)}_{y,x}\cdot h).
\end{array}$$
This shows that $\phi$ is a right $kG$-module homomorphism, and
hence a $kG$-Hopf bimodule isomorphism from $(kQ_1^c, G, r,
\overrightarrow \chi, u)$ to $^{\phi^{-1}}_{\phi}(kQ_1'{}^c, G', r',
\overrightarrow \chi', u')^{\phi^{-1}}_{\phi}$.

(ii) $\Rightarrow$ (i). Assume that there exist a Hopf algebra
isomorphism $\phi: kG\rightarrow kG'$ and a $kG$-Hopf bimodule
isomorphism $\psi: (kQ_1^c, G, r, \overrightarrow \chi,
u)\rightarrow\ ^{\phi^{-1}}_{\phi}(kQ_1'{}^c, G', r',
\overrightarrow \chi', u')^{\phi^{-1}}_{\phi}$. Then $\phi:
G\rightarrow G'$ is a group isomorphism. Let $C\in{\mathcal K}(G)$.
Then $\phi(u(C))$, $u'(\phi(C))\in\phi(C)\in{\mathcal K}(G')$, and
hence
$u'(\phi(C))=\phi(h_C)^{-1}\phi(u(C))\phi(h_C)=\phi(h^{-1}_Cu(C)h_C)$
for some $h_C\in G$. Since $\psi$ is a $kG'$-bicomodule isomorphism
from $^{\phi}(kQ_1^c, G, r, \overrightarrow \chi, u)^{\phi}$ to
$(kQ_1'{}^c, G', r', \overrightarrow \chi', u')$ and
$\phi(h^{-1}_Cu(C)h_C)=u'(\phi(C))$, by restriction one gets a
$k$-linear isomorphism
$$\psi_C:\ ^{h^{-1}_Cu(C)h_C}(kQ_1)^1\rightarrow\ ^{u'(\phi(C))}\! (kQ'_1)^1,\ x\mapsto \psi(x).$$
We also have a $k$-linear isomorphism
$$f_C:\ ^{u(C)}\! (kQ_1)^1\rightarrow\ ^{h^{-1}_Cu(C)h_C}(kQ_1)^1,\ x\mapsto h^{-1}_C\cdot x\cdot h_C.$$
Since $\phi(h^{-1}_Cu(C)h_C)=u'(\phi(C))$ and
$h^{-1}_CZ_{u(C)}h_C=Z_{h^{-1}_Cu(C)h_C}$, one gets
$\phi(h^{-1}_CZ_{u(C)}h_C)=Z_{u'(\phi(C))}$. Hence $\phi$ and $h_C$
induce an algebra isomorphism
$$\sigma_C:\ kZ_{u(C)}\rightarrow kZ_{u'(\phi(C))},\ h\mapsto \phi(h^{-1}_Chh_C).$$
Using the hypothesis that $\psi$ is a $kG$-bimodules homomorphism
from  $(kQ_1^c, G, r, \overrightarrow \chi, u)$   to
$_{\phi}(kQ_1'{}^c,$\ $ G', r', \overrightarrow \chi', u')_{\phi}$,
one can easily check that the composition $\psi_Cf_C$ is a right
$kZ_{u(C)}$-module isomorphism from $^{u(C)}\! (kQ_1)^1$ to
$(^{u'(\phi(C))}\! (kQ'_1)^1)_{\sigma_C}$. Since both  $^{u(C)}\!
(kQ_1)^1$   and $(^{u'(\phi(C))}\! (kQ'_1)^1)_{\sigma_C}$ \ \ are
pointed right $kZ_{u(C)}$-modules, they are semisimple
$kZ_{u(C)}$-modules for any $C\in{\mathcal K}_r(G)$. Moreover,
$ka^{(i)}_{u(C),1}$ and $ka^{(j)}_{u'(\phi(C)),1}$ are simple
submodules of $^{u(C)}\! (kQ_1)^1$ and $(^{u'(\phi(C))}\!
(kQ'_1)^1)_{\sigma_C}$, respectively, for any $i\in I_C(r)$ and
$j\in I_{\phi(C)}(r')$, where $C\in{\mathcal K}_r(G)$. Thus for any
$C\in{\mathcal K}_r(G)$, there exists a bijective map $\phi_C:
I_C(r)\rightarrow I_{\phi(C)}(r')$ such that $ka^{(i)}_{u(C),1}$ and
$(ka^{(\phi_C(i))}_{u'(\phi(C)),1})_{\sigma_C}$ are isomorphic right
$kZ_{u(C)}$-modules for any $i\in I_C(r)$, which implies
$\chi_C^{(i)}(h)=\chi'{}_{\phi(C)}^{(\phi_C(i))}(\phi(h^{-1}_Chh_C))$
for any $h\in Z_{u(C)}$ and $i\in I_C(r)$. It follows that ${\rm
RSC} (G, r , \overrightarrow{\chi}, u)$ $\cong $ ${\rm RSC} (G', r',
\overrightarrow{\chi'}, u')$. \ \ $\Box $

Up to now we have classified the {\rm PM} quiver Hopf algebras by
means of {\rm RSC}'s. In other words,  ramification systems with
characters uniquely determine  their corresponding PM quiver Hopf
algebras  up to graded Hopf algebra isomorphisms.

\begin {Example}\label{1.16} Assume that  $k$ is a field with char$(k)\not=2$.
Let $G=\{1, g\}\cong{\bf  Z}_2$ be the cyclic group of order $2$
with the generator $g$. Let $r$ be a ramification data of $G$ with
$r_1=m$ and $r_{\{g\}}=0$ and $(Q,G,r)$ be the corresponding Hopf
quiver, where $m$ is a positive integer. Then
$^1(Q_1)^1=\{a^{(i)}_{1,1}\mid i=1, 2, \cdots, m\}$,
$^g(Q_1)^g=\{a^{(i)}_{g,g}\mid i=1, 2, \cdots, m\}$, $^1(Q_1)^g$ and
$^g(Q_1)^1$ are two empty sets. For simplification, we write
$x_i=a^{(i)}_{1,1}$ and $y_i=a^{(i)}_{g,g}$ for any $1\leq i\leq m$.
Clearly, $Z_{u(\{1 \})}=G$ and $\widehat{G}=\{\chi_+, \chi_-\}$,
where $\chi_{\pm}(g)=\pm 1$. For any $0\leq n\leq m$, put
$\overrightarrow{\chi_n} \in (\hat G )^m $ with
$\chi_{n\{1\}}^{(i)}=\left\{\begin{array}{ll}
\chi_-,& \mbox{if } i>n;\\
\chi_+,& \mbox{otherwise}.
\end{array}\right.$
 Then $\{ {\rm RSC} (G, r, \overrightarrow{\chi _n}, u_n) \mid  n=0, 1, 2, \cdots, m \}$ are all non-isomorphic ${\rm RSC}$'s.
 Thus by Theorem \ref{4} we know that the path coalgebra
$kQ^c$ exactly admits $m+1$ distinct {\rm PM} co-path Hopf algebra
structures $kQ^c(G, r, \overrightarrow {\chi_n}, u_n)$, $0\leq n\leq
m$, up to graded Hopf algebra isomorphism. Now let $0\leq n\leq m$.
Then by Proposition \ref{1.10}, the $kG$-actions on $(kQ_1^c, G, r,
\overrightarrow {\chi_n}, u_n)$ are given by $g\cdot x_i=y_i,\ \
g\cdot y_i=x_i,\ \ 1\leq i\leq m;$$ \ \
$$x_i\cdot g=\left\{\begin{array}{rl}
-y_i,& \mbox{if }i>n,\\
y_i,& \mbox{otherwise},\\
\end{array}\right.\ \ \
 \ \ y_i\cdot g=\left\{\begin{array}{rl}
-x_i,& \mbox{if }i>n,\\
x_i,& \mbox{otherwise}.
\end{array}\right.$
Thus by \cite[p.245 or Theorem 3.8]{CR02}, the products of these
arrows $x_i, y_j$ in $kQ^c(G, r, \overrightarrow {\chi_n}, u_n)$ can
be described as follows. For any $i, j =1, 2, \cdots, m$, $x_i.x_j =
x_ix_j + x_jx_i,\ \ \  y_i.x_j=y_iy_j+y_jy_i,$$ \ \
$$x_i.y_j=\left\{\begin{array}{ccl}
-(y_iy_j + y_jy_i)&,& \mbox{if }i>n,\\
y_iy_j + y_jy_i&,& \mbox{otherwise},\\
\end{array}\right. \ \ \
\ y_i.y_j=\left\{\begin{array}{ccl}
-(x_ix_j + x_jx_i)&,& \mbox{if }i>n,\\
x_ix_j +x_jx_i&,& \mbox{otherwise},\\
\end{array}\right.$
where $x.y$ denotes the product of $x$ and $y$ in $kQ^c(G, r,
\overrightarrow {\chi_n}, u_n)$ for any $x, y\in kQ^c(G, r,
\overrightarrow {\chi_n}, u_n)$, $x_ix_j$ and $y_iy_j$ denote the
$2$-paths in the quiver $Q$ as usual for any $1\leq i, j\leq m$.
\end {Example}

\section{\bf  Multiple Taft algebras}\label{s2}

In this section  we discuss the {\rm PM} quiver Hopf algebras
determined by the {\rm RSC}'s with ${\mathcal K}_r(G)\subseteq
Z(G)$. We give the classification of {\rm PM} Nichols algebras and
multiple Taft algebras by means of  element system with characters
when $G$ is finite abelian group and $k$ is the complex field.

Let $r$ be a ramification data of $G$ and $(Q, G,r)$ be the
corresponding Hopf quiver. If $C$ contains only one element of $G$
for any $C\in {\mathcal K}_r(G)$, then $C=\{g\}$ for some $g\in
Z(G)$, the center of $G$. In this case, we say that the ramification
$r$ is {\it central}, and that ${\rm RSC} (G, r,
\overrightarrow{\chi}, u)$ {\it a central ramification system with
characters}, or a {\rm {\rm CRSC}} in short.  If ${\rm RSC} (G, r,
\overrightarrow{\chi}, u)$ is {\rm CRSC}, then the {\rm PM} co-path
Hopf algebra $kQ^c(G, r, \overrightarrow \chi, u)$ is called a {\it
multiple crown algebra} and $kG[kQ_1^c, G, r, \overrightarrow \chi,
u]$ is called a {\it multiple Taft algebra}.

\begin {Definition}\label{2.1}

$(G, \overrightarrow{g}, \overrightarrow{\chi}, J)$ is called an
{\it element system with characters} $($simply, {\rm ESC}$)$ if $G$
is a group, $J$ is a set, $\overrightarrow{ g } = \{g_i\} _{ i\in J}
\in Z(G)^J$ and $\overrightarrow{ \chi }= \{\chi_i\}_{ i \in J} \in
\widehat{G}^J $ with $ g_i \in Z(G)$ and $\chi _i \in \widehat G$.
${\rm ESC} (G, \overrightarrow{g}, \overrightarrow{\chi}, J)$ and
${\rm ESC} (G', \overrightarrow{g'}, \overrightarrow{\chi'}, J')$
are said to be isomorphic if there exist a group isomorphism $\phi:
G \rightarrow G'$ and a bijective map $\sigma: J\rightarrow J'$ such
that $\phi(g_i)=g'_{\sigma(i)}$ and $\chi'_{\sigma(i)}\phi=\chi_i$
for any $i \in J$.
\end {Definition}

$ {\rm ESC}(G,\overrightarrow{g},\overrightarrow{\chi},J)$ can be
written as ${\rm ESC}(G,g_i,\chi_i;i\in J)$ for convenience.
Throughout this paper, let $q_{ji}:=\chi_i(g_j)$, $q_i=q_{ii}$ and
$N_i$ be the order of $q_i$ ($N_i=\infty$ when $q_i$ is not a root
of unit, or $q_i=1$) for $i,j\in J$.

Let $H$ be a Hopf algebra with a bijective antipode $S$. A {\rm YD}
$H$-module $V$ is {\it pointed} if $V=0$ or $V$ is a direct sum of
one dimensional {\rm YD} $H$-modules. If $V$ is a pointed {\rm YD}
$H$-module, then the corresponding Nichols algebra ${\mathcal B}(V)$
is called a {\rm PM} {\it Nichols algebra}.

\begin{Lemma}\label{2.2}
Let $(V, \alpha ^-, \delta ^- )$ be a {\rm YD} $kG$-module. Then
$(V, \alpha ^-, \delta^- )$ is  a pointed  {\rm YD} $kG$-module  if
and only if $(V, \alpha ^-)$ is a pointed $($left$)$ $kG$-module
with $\delta^-(V)\subseteq kZ(G)\otimes V$.
\end{Lemma}

\begin{Lemma}\label{2.3} Assume that $G$ is a finite abelian group
of exponent $m$. If $k$ contains a primitive $m$-th root of 1, then
every {\rm YD} $kG$-module is  pointed  and  Nichols algebra  of
every {\rm YD} $kG$ -module is {\rm PM}.
\end{Lemma}
{\bf Proof.}  It follows from  Lemma \ref {1.2} and Lemma \ref
{2.2}. $\Box$

Let $(G, g_i, \chi_i; i\in J)$ be an ${\rm ESC}$. Let $V$ be a
$k$-vector space with ${\rm dim}(V)=|J|$. Let $\{x_i \mid i\in J\}$
be a basis of $V$ over $k$. Define a left $kG$-action and a left
$kG$-coaction on $V$ by
$$g\cdot x_i = \chi_i(g)x_i,\ \delta^-(x_i )= g_i \otimes
x_i,\ i\in J,\ g\in G.$$ Then it is easy to see that $V$ is a
pointed {\rm YD} $kG$-module and $kx_i$ is a one dimensional {\rm
YD} $kG$-submodule of $V$ for any $i\in J$. Denote by $V(G, g_i,
\chi_i; i\in J)$ the pointed {\rm YD} $kG$-module $V$. Note that
$V(G, g_i, \chi_i; i\in J)=0$ if $J$ is empty.

\begin{Proposition}\label{2.4}
$V$ is pointed {\rm YD} $kG$-module if and only if $V$ is isomorphic
to $V(G, g_i, \chi_i;$ $i\in J)$ for some  {\rm ESC} $(G, g_i,
\chi_i; i\in J)$.
\end{Proposition}

{\bf Proof.} If $V\cong V(G, g_i,\chi_i;i\in J)$ for some {\rm ESC}
$(G, g_i,\chi_i;i\in J)$ of $G$, then $V$ is obviously a pointed
{\rm YD} $kG$-module. Conversely, assume  that $V$ is a nonzero
pointed {\rm YD} $kG$-module. By Lemma \ref{2.2}, $V=\bigoplus_{g\in
Z(G)}V_g$ and $V_g=\{v\in V|\delta^-(v)= g\otimes v\}$ is a pointed
{\rm YD} $kG$-submodule of $V$ for any $g\in Z(G)$. Let $g\in Z(G)$
with $V_g\not=0$. Then $V_g$ is a nonzero pointed $kG$-module. Hence
there is a $k$-basis $\{x_i\mid i\in J_g\}$ such that $kx_i$ is a
$kG$-submodule of $V_g$ for any $i\in J_g$. It follows that there is
a character $\chi_i\in\widehat G$ for any $i\in J_g$ such that
$h\cdot x_i=\chi_i(h)x_i$ for all $h\in G$. For any $i\in J_g$, put
$g_i=g$. We may assume that these index sets $J_g$ are disjoint,
that is, $J_g\cap J_h=\emptyset$ for any $g\not= h$ in $Z(G)$ with
$V_g\not=0$ and $V_h\not=0$. Now let $J$ be the union of all the
$J_g$ with $g\in Z(G)$ and $V_g\not=0$. Then one can see that $(G,
g_i,\chi_i;i\in J)$ is an {\rm ESC} of $G$, and that $V$ is
isomorphic to $V(G, g_i,\chi_i;i\in J)$ as a {\rm YD} $kG$-module.\
\ $\Box$

Now we give the relation between ${\rm RSC}$ and ${\rm ESC}$. Assume
that $(G, g_i, \chi_i; i\in J)$ is an {\rm ESC} of $G$. We define a
binary relation $\sim$ on $J$ by
$$i \sim j \Leftrightarrow g_i = g_j,$$
where $i, j \in J$. Clearly, this is an equivalence relation. Denote
by $J/\!\!\sim$ the quotient set of $J$ modulo $\sim$. For any $i\in
J$, let $[i]$ denote the equivalence class containing $i$. That is,
$[i]:=\{j\in J \mid j\sim i \}$. Choose a subset $\bar J \subseteq
J$ such that the assignment $i\mapsto [i]$ is a bijective map from
$\bar J$ to $J/\!\!\sim$. That is, $\bar J$ is a set of
representative elements of the equivalence class.  Then
$J=\bigcup_{i\in \bar J}[i]$ is a disjoint union. Let $r_{ \{g_i\}}
=|[i]|$ for any $i\in\bar J$. Then $r=\sum_{i\in \bar J} r_{\{g_i\}}
\{g_i\}$ is a central ramification of $G$ with $I_{\{g_i\}}(r)=[i]$
for $i\in \bar J$. Moreover, ${\mathcal K}_r(G)=\{ \{ g_i \}\mid
i\in\bar J\}$. Put $\chi_{\{g_i\}}^{(j)}:=\chi_j$ for any $i\in\bar
J$ and $j\in [i]$. We obtain an ${\rm  CRSC}$, written ${\rm
CRSC}(G, r(g_i, \chi_i; i\in J),$ \ $ \overrightarrow \chi (g_i,
\chi_i; i\in J), u) $.
 Let $(Q, G, r)$ be the corresponding Hopf quiver
with $r=r(g_i,\chi_i;i\in J)$ and denote by $(kQ_1^c, g_i, \chi_i;
i\in J)$ the $kG$-Hopf bimodule $(kQ_1^c, G, r, \overrightarrow
\chi, u)$. Denote by $kQ^c(G, g_i, \chi_i; i\in J)$ and $kG[kQ_1^c,
G, g_i, \chi_i; i\in J]$ the corresponding multiple crown algebra
$kQ^c(G, r, \overrightarrow \chi, u)$ and multiple Taft algebra
$kG[kQ_1^c, G, r, \overrightarrow \chi, u]$, respectively. We also
denote by $kQ^a(G, g_i,\chi_i;i\in J)$, $kQ^{sc}(G, g_i,\chi_i;i\in
J)$, $kQ^s(G, g_i,\chi_i;i\in J)$ and    $(kG)^*[kQ_1^a, G,  g_i,
\chi_i; i\in J]$ the corresponding path Hopf algebra $kQ^a(G, r,
\overrightarrow \chi, u)$, semi-co-path Hopf algebra $kQ^{sc}(G, r,
\overrightarrow \chi, u)$, semi-path Hopf algebra $kQ^s(G, r,
\overrightarrow \chi, u)$ and one-type- path Hopf algebra
$(kG)^*[kQ_1^a, G, r, \overrightarrow \chi, u]$, respectively.

Conversely, assume that  $(G, r, \overrightarrow \chi, u) $ is a
${\rm CRSC}$. We may assume $I_{\{g\}}(r)\cap
I_{\{h\}}(r)=\emptyset$ for any $\{g\}\not=\{ h\}$ in ${\mathcal
K}_r(G)$. Let $J:=\bigcup_{\{g\}\in{\mathcal K}_r(G)}I_{\{g\}}(r)$.
For any $i\in J$, put $g_i:=g$ and $\chi_i:=\chi_{\{g\}}^{(i)}$ if
$i\in I_{\{g\}}(r)$ with $\{g\}\in{\mathcal K}_r(G)$. We obtain an
${\rm ESC}$, written ${\rm ESC} (G, \overrightarrow{g} (r,
\overrightarrow \chi, u),$ \ $ \overrightarrow{\chi } (r,
\overrightarrow \chi, u), J).$

From now on, assume $I_{\{g\}}(r)\cap I_{\{h\}}(r)= \emptyset $ for
any $\{g\}, \{h\}\in {\cal K}_r(G)$ with $\{g\} \not=\{ h\}$.  Note
that, in two cases above, for any $i, j \in I_{\{g\}}(r)$, we have
\begin {eqnarray} \label {e2.41} g_i = g_j =g, \ \  a^{(j)} _{g_i,
1} = a^{(j)}_{g_j,1}, \ \  \chi ^{(j)} _{\{g_i\}} = \chi
^{(j)}_{\{g_j\}} = \chi _j.
\end {eqnarray} Throughout this paper, let $E_j := a_{g_j,
1}^{(j)}$ for any $j\in J$.
\begin {Proposition} \label {2.5}    ${\rm CRSC}(G, r,
\overrightarrow \chi, u) $ \ $ \cong
 {\rm CRSC }(G', r',
\overrightarrow {\chi'}, u')$ \   if and only if \ ${\rm ESC} (G, $
\ $ \overrightarrow{g} (r, \overrightarrow \chi, u),$ \ $
\overrightarrow{\chi } (r, \overrightarrow \chi, u), J) $ \ $ \cong
$ \ $ {\rm ESC} (G', \overrightarrow{g'}(r', \overrightarrow {
\chi'}, $ \ $ u'),\overrightarrow{\chi ' } (r', \overrightarrow
{\chi'}, u'), J')$.
\end {Proposition}

{\bf Proof.} We use  notations above. If ${\rm CRSC}(G, r,
\overrightarrow \chi, u) \cong {\rm  CRSC} (G', r', \overrightarrow
\chi', u')$, then there exists a group isomorphism $\phi:
G\rightarrow G'$; for any $C\in{\mathcal K}_r(G)$, there exists a
bijective map $\phi_C : I_C(r) \rightarrow I_{\phi(C)}(r')$ such
that $\chi '{}_{\phi(C)}^{(\phi_C(i))}\phi =\chi_C^{(i)}$.
 Let $\sigma$ be the bijection from $J$ to $J'$ such
 that $\sigma (i) = \phi _C(i)$  for any $C \in {\cal K}_r(G)$, $i \in
 I_C(r)$. It is clear that
  $\phi(g_j)=g'_{\sigma(j)}$ and $\chi'_{\sigma(j)}\phi =\chi_j$ for any $j \in J.$  Thus ${\rm ESC} (G,
\overrightarrow{g} (r, \overrightarrow \chi, u),
\overrightarrow{\chi } (r, \overrightarrow \chi, u), J) $ \ $ \cong
$ \ $ {\rm ESC} (G', \overrightarrow{g'}(r', \overrightarrow {
\chi'}, u'), $ \ $\overrightarrow{\chi ' } (r', \overrightarrow
{\chi'}, u'), J')$.

Conversely, if ${\rm ESC} (G,$ \ $ \overrightarrow{g} (r,
\overrightarrow \chi, u),$ \ $ \overrightarrow{\chi } (r,
\overrightarrow \chi, u), J) $ \ \  $ \cong $ \ \ $ {\rm ESC} (G', $
$ \overrightarrow{g'}(r', $ $ \overrightarrow { \chi'}, u'), $ \ \ \
$\overrightarrow{\chi ' } (r', \overrightarrow {\chi'}, u'), J')$,
then there exist a group isomorphism $\phi: G \rightarrow G'$,  a
bijective map $\sigma: J\rightarrow J'$ such that
$\phi(g_j)=g'_{\sigma(j)}$ and $\chi'_{\sigma(j)}\phi =\chi_j $ for
any $j \in J.$ For any $C = \{g_i\} \in {\cal K}_r(G), j \in
I_C(r)$, we define $\phi _C = \sigma \mid _{I_C(r)}$  and  have
\begin {eqnarray*}\chi '{}_{\phi (C)}^{(\phi _C (j))} \phi &=& \chi '{}_{g'_{\sigma (i)}}^{(\sigma (j))}
\phi =\chi '{}_{\sigma (j)} \phi = \chi_ j = \chi ^{(j)}_C.\ \ \Box
\end {eqnarray*}

If  $ (G, g_i,\chi_i;i\in J)$ is  an {\rm ESC}, then $(kQ_1^c, G,
g_i, \chi_i; i\in J)$ is a $kG$-Hopf bimodule with module operations
$\alpha^-$ and $\alpha^+$. Define a new left $kG$-action on $kQ_1$
by
$$g\rhd x:=g\cdot x\cdot g^{-1},\ g\in G, x\in kQ_1,$$
where $g\cdot x=\alpha^-(g\otimes x)$ and $x\cdot
g=\alpha^+(x\otimes g)$ for any $g\in G$ and $x\in kQ_1$. With this
left $kG$-action and the original left (arrow) $kG$-coaction
$\delta^-$, $kQ_1$ is a YD $kG$-module. Let  $Q_1^1:=\{a\in Q_1 \mid
s(a)=1\}$. It is clear that $kQ_1^1$ is a YD $kG$-submodule of
$kQ_1$, denoted  by $(kQ_1^1, ad(G, g_i, \chi_i; i\in J))$.

\begin{Lemma}\label{2.6}
 $(kQ_1^1, ad(G, g_i,\chi_i;i\in J))$ and
$V(G, g_i,\chi^{-1}_i;i\in J)$   are isomorphic {\rm YD}
$kG$-modules.
\end{Lemma}

{\bf Proof.}   By definition, $V(G, g_i,\chi^{-1}_i;i\in J)$ has a
$k$-basis $\{x_i\mid i\in J\}$ such that $\delta^-(x_i)=g_i\otimes
x_i$ and $g\cdot x_i=\chi^{-1}_i(g)x_i$ for all $i\in J$ and $g\in
G$.  By Proposition \ref{1.10}, for any $j \in J, $ we have that
$g\rhd a_{g_j,1}^{(j)}=\chi_j(g^{-1})a_{g_j,1}^{(j)}$  and
$\delta^-(a_{g_j,1}^{(j)})=g_j\otimes a_{g_j,1}^{(j)}$. It follows
that there is a {\rm YD} $kG$-module isomorphism from $(kQ_1^1,
ad(G, g_i,\chi_i;i\in J))$ to $V(G, g_i,\chi^{-1}_i;i\in J)$ given
by $a_{g_j,1}^{(j)}\mapsto x_j$ for any  $j\in J$. $\Box$


\begin{Lemma}\label{2.7}
Let $B$ and $B'$ be two Hopf algebras with bijective antipodes. Let
$V$ be a {\rm YD} $B$-module. Assume that there is a Hopf algebra
isomorphism $\phi: B'\rightarrow B$. Then
$^{\phi^{-1}}_{\phi}{\mathcal B}(V)$ $\cong $ ${\mathcal
B}({}^{\phi^{-1}}_{\phi}V)$ as graded braided Hopf algebras in
$^{B'}_{B'} {\cal YD}.$
\end{Lemma}

\begin{Theorem}\label{5} Assume that
 $(G, g_i, \chi_i; i\in J)$ and $ (G', g_i', \chi_i'; i\in J')$ are
 two {\rm ESC}'s.
Then the following statements are equivalent:

{\rm(i)}    ${\rm ESC} (G, g_i, \chi_i; i\in J) \cong \ {\rm ESC}
(G', g_i', \chi_i'; i\in J')$.

 {\rm(ii)} \ ${ \rm CRSC}(G, r(G, g_i, \chi_i; i\in J),$ \ $ \overrightarrow \chi (G, g_i, \chi_i; i\in J), u) $  \ $ \cong $\ $
{\rm  CRSC }(G', r'(G', g_i', \chi_i'; i\in J'),$ \ $
\overrightarrow {\chi'}(G', g_i', \chi_i'; i\in J'), $\ $ u')$.

 {\rm(iii)} There is a Hopf
algebra isomorphism $\phi: kG\rightarrow kG'$ such that $V(G, g_i,
\chi_i;$ $i\in J)\cong\ ^{\phi^{-1}}_{\phi}V'(G' g_i', \chi_i';$
$i\in J')$ as {\rm YD} $kG$-modules.

 {\rm(iv)} There is
a Hopf algebra isomorphism $\phi: kG\rightarrow kG'$ such that
${\mathcal B}(V(G,g_i, \chi_i;$ $i\in J))\cong\
^{\phi^{-1}}_{\phi}{\mathcal B}(V'(G', g_i', \chi_i';$ $i\in J'))$
as graded braided Hopf algebras in $^{kG}_{kG}{\mathcal YD}$.

 {\rm(v)} There is a Hopf
algebra isomorphism $\phi: kG\rightarrow kG'$ such that $(kQ_1^1,
ad(G, g_i,\chi_i;i\in J))\cong\ ^{\phi^{-1}}_{\phi}(kQ'{}_1^1,
ad(G', g'_i,\chi'_i;i\in J'))$ as {\rm YD} $kG$-modules.

{\rm(vi)}    $kG[ kQ^c_1, G, g_i, \chi_i; i\in J] \cong \ kG'[
kQ'{}^c_1,G', g_i', \chi_i'; i\in J'].$

\end{Theorem}
{\bf Proof. }  We use he notations before Proposition \ref {2.5}.

(i) $\Rightarrow$ (ii). There exist a group isomorphism $\phi: G
\rightarrow G'$ and   a bijective map $\sigma: J\rightarrow J'$ such
that $\phi(g_j)=g'_{\sigma(j)}$ and $\chi'_{\sigma(j)}\phi =\chi_j $
for any $j \in J$. For any $C=\{g_i\} \in {\cal K}_r (G)$ and $j \in
I_{C}(r)$, we have $g_i = g_j$ and
\begin {eqnarray*} \chi_{ C}^{(j)} = \chi_{\{
g_i\}}^{(j)}&=& \chi_j= \chi'_{\sigma(j)}\phi\\
&=& {\chi'}_{\{g'_{\sigma (i)}\}}^{(\sigma (j))}\phi =
{\chi'}_{\{\phi(g _{i})\}}^{(\phi_{g_i} (j))}\phi =
{\chi'}_{\phi(C)}^{(\phi_{C} (j))}\phi \ .
\end {eqnarray*}

(ii) $\Rightarrow$ (i). There is  a group isomorphism $\phi: G
\rightarrow G'$ and   a bijection $\phi _C: \ I_C(r)\rightarrow I
_{\phi (C)} (r')$ such that $\chi '{}_{\phi (C)}^{\phi _C(j)} \phi =
\chi _C ^{(j)}$ for any $C = \{ g_i \} \in {\cal K}_r(G)$, $j \in
I_C(r)$. Define a map $\sigma : \ J \rightarrow J'$  such that
$\sigma \mid _{I_C(r)} = \phi _C$ for any $C \in {\cal K}_r(G)$.
Thus $\sigma $ is  bijective. For any $C= \{g_i\} \in {\cal K}_r(G)$
and $j \in I_{\{g_i\}}(r)$,  we have
$\phi(g_j)=\phi(g_i)=g'_{\sigma(j)}$  and
${\chi'}_{\sigma(j)}\phi={\chi'}_{\phi_{\{g_i\}}(j)}\phi$ $ =
{\chi'}_{\phi(\{g_i\})}^{(\phi _{\{g_i\}}(j))}\phi
=\chi_{\{g_{i}\}}^{(j)}=\chi_j$. This shows that ${\rm ESC} (G,
\overrightarrow{g}, $\ $\overrightarrow{\chi }, J)$ \ $\cong $ ${\rm
ESC} (G', \overrightarrow{g'}, $\ $\overrightarrow{\chi '}, J')$.

(i) $\Leftrightarrow$ (iii). Let  $V:= V(G, g_i, \chi_i;$ $i\in J)$
and $ V':= V'(G' g_i', \chi_i';$ $i\in J')$. By definition $V$ has a
$k$-basis $\{x_i\mid i\in J\}$ such that $g\cdot x_i=\chi_i(g)x_i$
and $\delta^-(x_i)=g_i\otimes x_i$ for any $i\in J$ and $g\in G$.
Similarly, $V'$ has a $k$-basis $\{y_j\mid j\in J'\}$ such that
$h\cdot y_j=\chi'_j(h)y_j$ and $\delta^-(y_j)=g'_j\otimes y_j$ for
all $j\in J'$ and $h\in G'$.

Assume ${\rm ESC} (G, \overrightarrow{g}, $\ $\overrightarrow{\chi
}, J)$ \ $\cong $ ${\rm ESC} (G', \overrightarrow{g'}, $\
$\overrightarrow{\chi '}, J')$. Then there is a group isomorphism
$\phi: G \rightarrow G'$ and a bijective map $\sigma : J \rightarrow
J'$ such that $\phi(g_i) = g'_{\sigma (i)}$ and
$\chi'_{\sigma(i)}\phi=\chi_i$ for any $i \in J$. Hence $\phi:
kG\rightarrow kG'$ is a Hopf algebra isomorphism. Define a
$k$-linear isomorphism $\psi: V\rightarrow V'$ by
$\psi(x_i)=y_{\sigma(i)}$ for any $i\in J$. Then it is
straightforward to check that $\psi$ is a {\rm YD} $kG$-module
homomorphism from $V$ to $^{\phi^{-1}}_{\phi}V'$.

Conversely, assume that $\phi: kG\rightarrow kG'$ is a Hopf algebra
isomorphism and $\psi: V\rightarrow\ ^{\phi^{-1}}_{\phi}V'$ is a
{\rm  YD} $kG$-module isomorphism. Then $\phi: G\rightarrow G'$ is a
group isomorphism. We use the notations in the proof of Lemma
\ref{2.2} and the notations above. Then $\psi(V_g)=V'_{\phi(g)}$ for
any $g\in G$. Since $V=\bigoplus_{i\in\bar J}V_{g_i}$ and
$V'=\bigoplus_{i\in\bar{J'}}V_{g'_i}$, there is a bijection $\tau:
\bar J \rightarrow \bar J'$ such that
$\psi(V_{g_i})=V'_{g'_{\tau(i)}}$ for any $i\in\bar J$. This shows
that $\phi(g_i)=g'_{\tau(i)}$ and $V_{g_i}\cong\
_{\phi}(V'_{g'_{\tau(i)}})$ as left $kG$-modules for any $i\in \bar
J$. However, $V_{g_i}$ is a pointed $kG$-module and $kx_j$ is its
one dimensional submodule for any $j\in[i]$. Similarly,
$V_{g'_{\tau(i)}}$ is a pointed $kG'$-module and $ky_j$ is its one
dimensional submodule for any $j\in[\tau(i)]$. Hence there is a
bijection $\phi_{\{g_i\}}: [i]\rightarrow [\tau(i)]$ for any $i\in
\bar J$ such that $kx_j$ and $_{\phi}(ky_{\phi_{\{g_i\}}(j)})$ are
isomorphic $kG$-modules for all $j\in[i]$. This implies that
$\chi_j=\chi'_{\phi_{\{g_i\}}(j)}\phi$ for all $j\in[i]$ and
$i\in\bar J$. Then the same argument as in the proof of (ii)
$\Rightarrow$ (i) shows that ${\rm ESC} (G, \overrightarrow{g}, $\
$\overrightarrow{\chi }, J)$ \ $\cong $ ${\rm ESC} (G',
\overrightarrow{g'}, $\ $\overrightarrow{\chi '}, J')$. (iii)
$\Leftrightarrow$ (iv)  It follows from Lemma \ref{2.7}. \ (iii)
$\Leftrightarrow$ (v) It follows from Lemma \ref {2.6}. \ (ii)
$\Leftrightarrow$ (vi) It follows from Theorem \ref {4}. $\Box$

Up to now we have classified  Nichols algebras and {\rm YD} modules
over finite abelian group and  the complex field up to isomorphisms,
which are under  means of Theorem  \ref {5} (iv)(iii), respectively.
In fact, we can explain these facts above by introducing some new
concepts about isomorphisms. For convenience, if $B$ is a Hopf
algebra and $M$ is a $B$-Hopf bimodule, then  we say that $(B, M)$
is a  Hopf bimodules. For any two Hopf bimodules $(B,M)$ and $(B',
M')$, if $\phi$ is a Hopf algebra homomorphism from $B$ to $B'$
 and $\psi$ is simultaneously a $B$-bimodule homomorphism from
$M$ to $_\phi M'{}_\phi $ and a $B'$-bicomodule homomorphism from
$^\phi M ^\phi $ to $M'$, then $(\phi, \psi)$ is called a pull-push
Hopf bimodule homomorphism. Similarly, we say that $(B, M)$ and $(B,
X)$ are a {\rm YD} module and a {\rm YD} Hopf algebra if $M$ is a
{\rm YD} $B$-module and $X$ is a braided Hopf algebra in
Yetter-Drinfeld category $^B_B {\cal YD}$, respectively.
 For any two
{\rm YD} modules $(B,M)$ and $(B', M')$, if $\phi$ is a Hopf algebra
homomorphism from $B$ to $B'$,  and $\psi$ is simultaneously a left
$B$-module homomorphism from  $M$ to $_\phi M' $ and a left
$B'$-comodule homomorphism from $^\phi M  $ to $M'$, then $(\phi,
\psi)$ is called a pull-push {\rm YD} module homomorphism. For any
two {\rm YD} Hopf algebra $(B,X)$ and $(B', X')$, if $\phi$ is a
Hopf algebra homomorphism from $B$ to $B'$,   $\psi$ is
simultaneously a left $B$-module homomorphism from  $X$ to $_\phi X'
$ and a left $B'$-comodule homomorphism from $^\phi X  $ to $X'$,
meantime, $\psi$ also is algebra and coalgebra homomorphism from $X$
to $X'$, then $(\phi, \psi)$ is called a pull-push {\rm YD} Hopf
algebra homomorphism.

Consequently, we have classified  Nichols algebras  over finite
abelian group and the complex field up to  pull-push graded {\rm YD}
Hopf algebra isomorphisms and {\rm YD} modules  over finite abelian
group and  the complex field up to pull-push  {\rm YD} module
isomorphisms, respectively. In other words,  element systems with
characters uniquely determine  their corresponding Nichols algebras
and {\rm YD} modules up to their isomorphisms.

\section    {The relation between quiver Hopf algebras and quotients  of free algebras
}\label {s3}

In this section we show that the diagram of a quantum weakly
commutative multiple Taft algebra is   not only a Nichols algebra
but also a quantum linear space in $^{kG}_{kG}{\cal YD}$; the
diagram of a semi-path Hopf algebra of ${\rm ESC}$   is  a quantum
tensor algebra in $^{kG}_{kG}{\cal YD}$; the quantum enveloping
algebra of a complex semisimple Lie algebra is a quotient of a
semi-path Hopf algebra.

\subsection {The structure of multiple Taft algebras and semi-path Hopf algebras}\label {s3.1}
Assume that $H=\bigoplus_{i\geq 0}H_{(i)}$ is a graded Hopf algebra
with invertible antipode $S$. Let $B=H_{(0)}$, and let $\pi_0:
H\rightarrow H_{(0)}=B$ and $\iota_0: B=H_{(0)}\rightarrow H$ denote
the canonical projection and injection. Set
$\omega:=id_H*(\iota_0\pi_0S): H\rightarrow H$. Then it is clear
that $(H, \delta ^+, \alpha ^+)$ is a right $B$-Hopf module with
$\delta^+:=({\rm id}\otimes\pi_0)\Delta $ and $\alpha^+:=\mu({\rm
id}\otimes\iota_0)$.   Let $R:=H^{co B}:=\{h\in
H\mid\delta^+(h)=h\otimes 1\}$, which is a graded subspace of $H$.
Then it is known that $R={\rm Im}(\omega)$ and $\Delta(R)\subseteq
H\otimes R$. Hence $R$ is a left coideal subalgebra of $H$, and so
$R$ is a left $H$-comodule algebra. It is well known that $R$ is a
graded braided Hopf algebra in $^B_B{\mathcal YD}$ with the same
multiplication, unit and counit as in $H$,  the comultiplication
$\Delta_R=(\omega \otimes id )\Delta$ , where the left $B$-action
$\alpha _R$ and left $B$-coaction $\delta _R$ on $R$ are given by
\begin {eqnarray}\label {e3.11}
\alpha _R (b \otimes x) = b\rightharpoonup_{ad}x=\sum
b_{(1)}xS(b_{(2)}),\ \delta_R^-(x)=\sum \pi_0(x_{(1)})\otimes
x_{(2)},\ b\in B,\ x\in R
\end {eqnarray}( see the proof of \cite [Theorem 3]{Ra85}). $R$ is
called the diagram of $H$, written $diag (H)$. Note that diagram $R$
of $H$ is dependent on  the gradation of $H$.  By \cite [Theorem
1]{Ra85}, the biproduct of $R$ and $B$ is a Hopf algebra, written $R
^{\delta _R}_{\alpha _R} \# B$, or $R \# B$ in short. The biproduct
$R ^{\delta _R}_{\alpha _R} \# B$ is also called the bosonization of
$R$. Furthermore, we have the following well known result.

\begin{Theorem}\label{6}
$($see \cite[p.1530]{Ni78}, \cite{AS98a} and \cite{Ra85}$)$  Under
notations  above, if $H=\bigoplus_{i \geq 0 }H _{(i)}$ is a graded
Hopf algebra, then $R $ is a graded braided Hopf algebra in
$^B_B{\mathcal YD}$ and $ R\#B\cong H$ as graded Hopf algebras,
where the isomorphism is $\alpha ^+  := \mu _H(id _H \otimes \iota
_0) $.
\end{Theorem}

{ \bf Remark}: If  $A$ be a Hopf algebra whose coradical $A_{0}$ is
a Hopf subalgebra, then it is clear that $H:=gr A$ is a graded Hopf
algebra. The diagram of $H$ with respect to gradation of $gr A$ is
called the diagram  of $A$ in \cite [Introduction ]{AS98b}.

\begin {Lemma}\label {3.1} (i) Assume that  $H$ and $H'$ are two graded Hopf
algebras with $B= H_{(0)}$ and $B' = H'_{(0)}.$  Then $H \cong H'$
as graded Hopf algebras if and only if there exists a Hopf algebra
isomorphism $\phi: B\rightarrow B'$ such that $ diag (H)\cong\
_{\phi} ^{\phi ^{-1}}diag (H')$ as {\rm YD} $B$-modules and as
graded  braided Hopf algebras in $^B_B {\cal YD}$.

(ii) Let $B$ and $B'$ be two  Hopf algebras. Let $M$ and $M'$ be
$B$-Hopf bimodule and  $B'$-Hopf bimodule, respectively. Then $B[M]
\cong B'[M']$ as graded Hopf algebras if and only if there exists a
Hopf algebra isomorphism $\phi: B\rightarrow B'$ such that $ diag
(B[M])\cong\ _{\phi} ^{\phi ^{-1}}diag (B'[M'])$ as {\rm YD}
$B$-modules and as graded braided Hopf algebras in $^B_B {\cal YD}$.

\end {Lemma}
{\bf Proof.} (i) Assume that $\xi$ is a graded Hopf algebra
isomorphism from $H$ to $H'$. Let $R := diag (H)$, $R' := diag
(H')$, $\phi := \xi \mid _ B$ and $\psi := \xi \mid _ R$. It is easy
to check that $\psi$ is the map required.

 Conversely, by Theorem
\ref {6}, $R \# B\cong H$ and $R'\# B' \cong H'$ as graded Hopf
algebras.  Let $\xi $ be a linear map from $ R \# B$ to $R' \# B'$
by sending $r \# b$ to $\psi (r) \# \phi (b)$ for any $r \in R$,
$b\in B$. Let $\nu $ be a linear map from $ R' \# B'$ to $R \# B$ by
sending $r' \# b'$ to $\psi ^{-1} (r') \# \phi ^{-1} (b')$ for any
$r' \in R'$, $b'\in B'$. Obviously, $\nu$ is the inverse of $\xi.$
Since $\psi$ is graded, so is $\xi$.

 Now we show that $\xi$ is an algebra homomorphism.
For any $r, r'\in R, b, b' \in B$, see
\begin {eqnarray*}
\xi \mu _{R\# B}( (r \# b) \otimes (r' \# b'))&=&  \psi (r (b _{(1)}
\cdot r'))\# \phi (b_{(2)}b')\\
&=&\psi (r) (\phi (b _{(1)})\cdot \psi (r' ) ) \# \phi (b_{(2)})
\phi (b) \ \\
&& (\hbox {since } \psi  \hbox { is a pullback module homomorphism }
\\
&& \hbox {
and  an algebra homomorphism. } ) \\
&=&\mu _{R'\# B'} (\xi (r \# b) \otimes \xi (r' \# b')).
\end {eqnarray*}
Similarly, we can show that $\xi$ is a coalgebra homomorphism.

(ii) It follows from (i). $\Box$



\begin{Lemma}\label{3.2} (i) $(kQ^c (G, r, \overrightarrow \chi, u))^{co \ kG} = span \{ \beta \mid \beta \hbox { is a path with  }
s(\beta ) =1 \}$. (ii) $(kG[kQ_1^c, $ $G, r, \overrightarrow \chi, u
])^{co \ kG}$ is the subalgebra of $kG[kQ^c_1, G, r, \overrightarrow
\chi, u]$ generated by $Q_1^1$ as algebras. (iii)  $(kG[kQ_1^c, $ $
G, r, \overrightarrow \chi, u] )^{co \ kG}$ $ \# kG \cong kG[kQ_1^c,
G, r, \overrightarrow \chi, u ]$ as graded Hopf algebra isomorphism.
(iv) $(kQ^s(G, r, \overrightarrow \chi, u))^{co \ kG}$ is the
subalgebra of $kG^s( G, r, \overrightarrow \chi, u)$ generated by
$Q_1^1$ as algebras.

\end{Lemma}

{\bf Proof.} (i)
 For a path $\beta $, see that
 \begin {eqnarray*}
 \delta ^+ (\beta ) = (id \otimes \pi _0)\Delta (\beta ) &=&  \beta \otimes s(\beta). \end {eqnarray*} This implies
$(kQ^c)^{co \ kG} = span \{ \beta \mid \beta \hbox { is a path with
} s(\beta ) =1 \}$.

(ii)  Since every path generated by arrows  in $Q_1^1$ is of start
vertex 1, this  path belongs to $(kG[kQ_1^c,G, r, \overrightarrow
\chi, u ])^{co \ kG}$. Let $R := (kG[kQ_1^c, G, r, \overrightarrow
\chi, u])^{co \ kG}$ and
 $A: = $ the subalgebra of
$kG[kQ^c_1, G, r, \overrightarrow \chi, u]$ generated by $Q_1^1$ as
algebras. Obviously, $A
 \subseteq
R$. It
 is clear that $ \alpha ^+ (R \# kG)  = \alpha ^+ (A \# kG) = kG[kQ_1^c,G, r, \overrightarrow \chi, u]$ and $\alpha ^+$ is injective. Thus $R \# kG  = A \#
 kG$ and $R = A$.

(iii) It follows from   Theorem \ref{6}.

(iv) We first show   $(kQ^s)^{co\ kG} =$  $span \{ \beta \mid \beta
=1,  \hbox {or }  \beta = \beta _n \otimes _{kG} \beta _{n-1}
\otimes _{kG} \cdots \otimes _{kG} \beta _1 $ \ $ \hbox {with }  $ $
\prod _{i =1}^n s( \beta _i ) = 1 \hbox { and } \beta _i \in Q_1, \
i = 1, 2, \cdots, n; n \in {\mathbb Z}^+ \}$. Indeed, obviously
right hand side $ \subseteq $ the left hand side. For any $\beta =
\beta _n \otimes _{kG} \beta _{n-1} \otimes _{kG} \cdots \otimes
_{kG} \beta _1 $\ $ \hbox {with } \beta _i \in Q_1 $, called a
monomial, define $s(\beta ) = \prod _{i =1}^n s( \beta _i )$. For
any $0 \not= u \in (kQ^s)^{co\ kG}$ with $u \not\in kG$, there exist
linearly independent monomials $u_1, u_2, \cdots, u _n$ such that
 $u = \sum
_{i=1}^n b_iu_i$ with $0\not= b_i \in k$  for $i = 1, 2, \cdots n$.
See $\delta ^+ (u)$ $= \sum _{i=1}^n b_iu_i \otimes s(u_i) $ $= u
\otimes 1$. Consequently, $s(u_i) =1$ for $i =1,2, \cdots, n$. This
implies that $u $ belongs to the right hand side.

For any $\beta = \beta _n \otimes _{kG} \beta _{n-1} \otimes _{kG}
\cdots \otimes _{kG} \beta _1 \hbox {with }   \prod _{i =1}^n s(
\beta _i ) = 1 \hbox { and } \beta _i \in Q_1, \ i = 1, 2, \cdots,
n, $ we show that $\beta $ can be written as multiplication of
arrows in $Q_1^1$ by induction. When $n=1$, it is clear. For $n >1$,
see $\beta = \beta _n \otimes _{kG} \beta _{n-1} \otimes _{kG}
\cdots \otimes _{kG} \beta _2 \cdot  s(\beta _1)\otimes _{kG} (
s(\beta _1)^{-1}\cdot \beta _1 )$. Thus  $\beta $ can be written as
multiplication of arrows in $Q_1^1$. Consequently, we complete the
proof of (iv). $\Box$

Recall that a braided algebra $A$ in braided tensor category $({\cal
C}, C)$ with braiding $C$ is said to be  braided commutative or
quantum commutative, if $ab = \mu C( a \otimes b)$ for any $a, b\in
A$. An $ {\rm ESC}(G, g_i, \chi _i; i\in J ) $ is said to be quantum
commutative if
$$\chi _i (g_j) \chi _j (g_i)=1$$ for any $i, j \in J$. An $ {\rm ESC}(G, g_i, \chi _i; i\in J ) $ is said to be quantum weakly
commutative if
$$\chi _i (g_j) \chi _j (g_i)=1$$ for any $i, j \in J$ with
$i\not=j$.

\begin {Lemma} \label {3.3}

(i) $ {\rm ESC}(G, g_i, \chi _i; i\in J ) $ is   quantum weakly
commutative if and only if in $diag (kG [kQ_1^c, G, g_i, \chi _i; i
\in J])$,
\begin {eqnarray}\label {3.10''e1} E_i\cdot E_j = \chi _j
(g_i^{-1}) E_j\cdot E_i \end {eqnarray} for any $i, j \in J$ with $i
\not= j.$

(ii) $diag (kG [kQ_1^c, G, g_i, \chi _i; i \in J])$ is quantum
commutative in $^{kG}_{kG} {\cal YD}$ if and only if $ {\rm ESC}(G,
g_i, \chi _i; i\in J ) $ is   quantum commutative.

\end {Lemma}

For any positive integers $m$ and $n$, let
$$D_n^{n+m}=\{d=(d_{n+m},d_{n+m-1},\cdots ,d_1)\mid d_i=0\mbox{ or
}1, \sum_{i=1}^{n+m}d_i=n\}.$$ Let $d\in D_n^{n+m}$ and let
$A=a_na_{n-1}\cdots a_1\in Q_n$ be an $n$-path. We define a sequence
$dA=((dA)_{n+m},\cdots,(dA)_1)$ by
$$(dA)_i=\left\{\begin{array}{ccl}
t(a_{d(i)})&,& \mbox{ if }d_i=0;\\
a_{d(i)}&,& \mbox{ if }d_i=1,\\
\end{array}\right.$$
where $1\leq i\leq n+m$ and $d(i)=\sum_{j=1}^id_j$. Such a sequence
$dA$ is called an $n+m$-thin splits of the $n$-path $A$. Note that
if $d(i)=0$ then we regard $t(a_{d(i)})=s(a_1)$, since
$s(a_{d(i)+1})=s(a_1)$ in this case.


 If $0 \not= q \in
k$ and $0 \le i \le n < ord (q) $ (the order of $q$), we set
$(0)_{q}! =1$,
$$ \left ( \begin {array} {c} n\\
i
\end {array} \right )_q
 =  \frac{(n)_{q}!}{(i)_{q}!(n - i)_{q}!}, \quad \hbox {where
}(n)_{q}! = \prod_{1 \le i \le n} (i)_{q}, \quad (n)_{q} =
\frac{q^{n} - 1}{q - 1}.$$ In particular, $(n)_q = n$ when $q=1.$

\begin {Lemma}\label{3.4} In $kQ^c(G, r, \overrightarrow \chi, u)$, assume $\{ g\} \in
{\mathcal K}_r (G)$ and  $j\in I_{\{g\}}(r)$. Let $q:= \chi
^{(j)}_{\{g\}} (g)$.
 If $i_1, i_2, \cdots, i_m$ be non-negative integers, then
$$\begin{array}{rcl}
a^{(j)}_{g^{i_m +1},g^{i_m }}\cdot a^{(j)}_{g^{i_{m-1}
+1},g^{i_{m-1}}}\cdot\cdots\cdot a^{(j)}_{g^{i_1 +1},g^{i_1 }}
&=&q^{\beta_m}(m)_q! P^{(j)}_{g^{\alpha_m}}(g,m)\\
\end{array}$$
where $\alpha _m = i_1 + i_2 + \cdots + i_m $, $P^{(j)}_h(g,m) =$ \
$ a^{(j)}_{g^mh,g^{m-1}h}a^{(j)}_{g^{m-1}h, g^{m-2}h}\cdots
a^{(j)}_{gh,h}$, $\beta_1=0$ and $\beta_m=\sum_{j
=1}^{m-1}(i_1+i_2+\cdots+i_j )$ if $m>1$.
\end{Lemma}
{\bf Proof.} We prove the  equality by induction on $m$. For $m=1$,
it is easy to see that the equality holds. Now suppose $m>1$. We
have
$$\begin{array}{rl}
&a^{(j)}_{g^{i_m +1},g^{i_m }}\cdot a^{(j)}_{g^{i_{m-1}
+1},g^{i_{m-1}}}\cdot\cdots\cdot a^{(j)}_{g^{i_1 +1},g^{i_1 }}\\
=&a^{(j)}_{g^{i_m +1},g^{i_m }}\cdot(a^{(j)}_{g^{i_{m-1}
+1},g^{i_{m-1}}}\cdot\cdots\cdot a^{(j)}_{g^{i_1 +1},g^{i_1 }})\\
=&q^{\beta_{m-1}}(m-1)_q! a^{(j)}_{g^{i_m +1},g^{i_m }}\cdot
P^{(j)}_{g^{\alpha_{m-1}}}(g,m-1) \ \ \ ( \hbox {by inductive
assumption })\\
=&q^{\beta_{m-1}}(m-1)_q! a^{(j)}_{g^{i_m +1},g^{i_m}}\cdot
(a^{(j)}_{g^{\alpha_{m-1}+m-1},g^{\alpha_{m-1}+m-2}}\cdots
a^{(j)}_{g^{\alpha_{m-1}+1},g^{\alpha_{m-1}}}) \\
 =&q^{\beta_{m-1}}(m-1)_q! \sum_{l=1}^m[(g^{i_m+1}\cdot
a^{(j)}_{g^{\alpha_{m-1}+m-1},g^{\alpha_{m-1}+m-2}})
\cdots(g^{i_m+1}\cdot a^{(j)}_{g^{\alpha_{m-1}+l},g^{\alpha_{m-1}+l-1}})\\
&(a^{(j)}_{g^{i_m+1},g^{i_m}}\cdot g^{\alpha_{m-1}+l-1})
(g^{i_m}\cdot
a^{(j)}_{g^{\alpha_{m-1}+l-1},g^{\alpha_{m-1}+l-2}})\cdots
(g^{i_m}\cdot a^{(j)}_{g^{\alpha_{m-1}+1},g^{\alpha_{m-1}}})]\\
&\ \ \ ( \hbox {by \cite[Theorem 3.8]{CR02} })
\\
=&q^{\beta_{m-1}}(m-1)_q!
\sum_{l=1}^m[a^{(j)}_{g^{\alpha_m+m},g^{\alpha_m+m-1}}
\cdots a^{(j)}_{g^{\alpha_m+l+1},g^{\alpha_m+l}}\\
&(\chi_{\{g\}}^{(j)}(g^{\alpha_{m-1}+l-1})a^{(j)}_{g^{\alpha_m+l},g^{\alpha_m+l-1}})
a^{(j)}_{g^{\alpha_m+l-1},g^{\alpha_m+l-2}}\cdots
a^{(j)}_{g^{\alpha_m+1},g^{\alpha_m}}]    \ \ ( \hbox {by } Proposition  \ref {1.10}) \\
=&q^{\beta_{m-1}}(m-1)_q! \sum_{l=1}^m
q^{\alpha_{m-1}+l-1}P^{(j)}_{g^{\alpha_m }}(g,m)\\
=&q^{\beta_{m-1}+\alpha_{m-1}}(m)_q!P^{(j)}_{g^{\alpha_m }}(g,m)\\
=&q^{\beta_m}(m)_q! P^{(j)}_{g^{\alpha_m }}(g,m). \ \ \Box
\end{array}$$

\begin {Lemma} \label {3.5} (See \cite [Lemma 3.3]{AS98b}) Let $B$ be a Hopf algebra and $R$  a braided Hopf algebra in
${}_B^B {\cal YD}$ with a linearly independent set $ \{ x_{1} \dots,
x_{t} \}$ $\subseteq $ $P(R)$.  Assume that there exist $g_{j} \in
G(B)$ (the set of all group-like elements in $B$) and  $\chi_{j} \in
Alg(B,k)$ such that
 $$\delta (x_{j}) = g_{j}\otimes x_{j}, \ h\cdot  x_{j} = \chi_{j}(h)x_{j},
 \hbox { for all } h\in B, j =1, 2, \cdots, t . $$
Then \begin{eqnarray*}\{x_1^{m_1} x_2^{m_2}\cdots x_t^{m_t} \mid 0
\leq m_j<N_j, 1\le j \le t \}. \end{eqnarray*} is linearly
independent, where $N_i$ is the order of $q_i := \chi _i (g_i)$  \ (
$N_i = \infty $ when $q_i$ is not a root of unit, or $q_i =1$ ) for
$1 \le i \le t.$
\end {Lemma}

 Assume that $(G, g_i, \chi_i; j\in J)$ is an ${\rm ESC}.$ Let
${\cal T} (G, g_i, \chi_i; j\in J)$ be the free algebra generated by
set $ \{ x_j \mid j\in J \}$. Let ${\cal S} (G, g_i, \chi_i; j\in
J)$ be the algebra generated by set $ \{ x_j \mid j\in J \}$ with
relations
\begin {eqnarray} \label {qse1} x_{i}x_{j} = \chi_{j}(g_{i})
x_{j}x_{i}  \ \ \ \hbox { for any } i, j \in J \hbox { with } i
\not= j.
\end {eqnarray}
Let ${\cal R} (G, g_i, \chi_i; j\in J)$ be the algebra generated by
set $ \{ x_j \mid j\in J \}$  with relations
\begin {eqnarray} \label {qlse1} x_l ^{N_l}=0,\  x_{i}x_{j} = \chi_{j}(g_{i})
x_{j}x_{i}  \ \ \ \hbox { for any } i, j,l  \in J \hbox { with } N_l
< \infty,  i \not= j.
\end {eqnarray}
 Define their coalgebra operations and $kG$-(co-)module operations  as follows:
\begin {eqnarray} \label {e3.522}\Delta x_i = x_i \otimes 1 + 1 \otimes x_i, \ \ \epsilon (x_i)
=0, \ \ \delta ^-(x_{i}) = g_{i}\otimes x_{i}, \qquad h \cdot x_{i}
= \chi_{i}(h)x_{i}.\end {eqnarray} $ {\cal T} (G, g_i, \chi_i; j\in
J)$ is called a quantum tensor algebra in $^{kG}_{kG} {\cal YD} $,
${\cal S} (G, g_i, \chi_i; j\in J)$ is called a quantum symmetric
algebra in $^{kG}_{kG} {\cal YD} $ and ${\cal R} (G, g_i, \chi_i;
j\in J)$ is called a quantum linear space in $^{kG}_{kG} {\cal YD}
$. Note that when ${\rm ESC} (G, g_i, \chi _i; i \in J)$ is quantum
weakly commutative with finite $J$ and  finite $N_j$ for any $j\in
J$, the definition of quantum linear space is the same as in \cite
[Lemma 3.4]{AS98b}. Obviously, if $N_i$ is infinite for all $i \in
J$, then ${\cal S} (G, g_i, \chi_i; j\in J)$ = ${\cal R} (G, g_i,
\chi_i; j\in J)$.

\begin{Theorem}\label{7}
Assume that ${\rm ESC} (G, g_i, \chi _i; i \in J)$ is quantum weakly
commutative.  Let $\prec$ be  a total order of $J$. Then

{\rm(i)} \ The multiple Taft algebra $kG[kQ_1^c,  G,  g_i, \chi_i;
i\in J]$ has a $k$-basis
\begin{eqnarray*}\label{basis} \begin {array} {c} \{  g\cdot E_{\nu _1}^{m_1}\cdot
 E_{\nu _2}^{m_2}\cdot\cdots\cdot  E_{\nu _t}^{m_t} \mid 0 \le
m_j < N_j; \nu _j \prec \nu _{j+1}, \nu _j \in J,  j= 1, 2, \cdots ,
t;  t \in {\mathbb Z}^+, g\in G\}. {} \end {array} \end{eqnarray*}
Moreover, $kG[kQ_1^c, G g_i, \chi_i; i\in J]$ is finite dimensional
if and only if $\mid \! G \! \mid,$  $\mid \! J \! \mid $  and $N_j
$ are finite for any $j\in J$. In this case, ${\rm dim}_k(kG[kQ_1^c,
G, g_i, \chi_i; i\in J])=|G|N_1N_2\cdots N_t$ with $J = \{1, 2,
\cdots, t\}$.

 {\rm(ii)}\ $diag ( kG[kQ_1^c, G, g_i, \chi _i; i\in J])$ has a $k$-basis
\begin{eqnarray}\label{basis2}\{ E_{\nu _1}^{m_1}\cdot
 E_{\nu _2}^{m_2}\cdot\cdots\cdot  E_{\nu _t}^{m_t} \mid 0 \le
m_j < N_j; \nu _j \prec \nu _{j+1}, \nu _j \in J,  j= 1, 2, \cdots ,
t;  t \in {\mathbb Z}^+\}.\end{eqnarray}

{\rm(iii)} \ $diag ( kG[kQ_1^c, G, g_i, \chi _i; i\in J])$ is a
Nichols algebra in   ${}_{kG}^{kG}{\cal YD}$  and   ${\cal R} (G,
g_i, \chi_i ^{-1}; j\in J) \cong diag ( kG[kQ_1^c, G,  g_i, \chi _i;
i\in J])$ as graded braided Hopf algebras in ${}_{kG}^{kG}{\cal
YD}$,  by sending $x_j$ to $a_{g_j, 1}^{(j)}$ for any $j\in J$.

{\rm(iv)} \ ${\cal T} (G, g_i, \chi_i ^{-1}; j\in J) \cong  diag (
kQ^s( g_i, \chi _i; i\in J))$ as graded braided Hopf algebras in
${}_{kG}^{kG}{\cal YD}$ algebras, by sending $x_j$ to $a_{g_j,
1}^{(j)}$ for any $j\in J$.

{\rm(v)} \  $kQ^s ( G, g_i, \chi_i; i\in J) \cong {\cal T} (G, g_i,
\chi_i ^{-1}; j\in J)\# kG$ \ as graded Hopf algebras and  $kQ^s (G,
g_i, \chi_i; i\in J)$ has a $k$-basis
\begin{eqnarray*} \{  g\cdot E_{\nu _1}\otimes _{kG}
 E_{\nu _2}\otimes _{kG}\cdots \otimes _{kG} E_{\nu _t} \mid  \nu _j \in J,  j= 1, 2, \cdots
, t;  t \in {\mathbb Z}^+ \cup \{0\}, g\in G\}, \end{eqnarray*}
where $g\cdot E_{\nu _1}\otimes _{kG}
 E_{\nu _2}\otimes _{kG}\cdots \otimes _{kG} E_{\nu _t}=g$ when $t=0.$

 Note that (iv) and (v) still hold without quantum weakly commutative
condition.
 \end{Theorem}
{\bf Proof.} (ii)    Since $(N_j)_{q_j}!=0$, it follows from Lemma
\ref{3.4} that $E_j^{N_j}=0$ when $N_j <\infty$. By Lemma \ref
{3.3}, $E_i\cdot E_j = \chi _j (g_i ^{-1}) E_j\cdot E_i$ for any $i,
j \in J$ with $i \not= j$. Considering Lemma \ref{3.5}, we complete
the proof.

(iii) By Lemma \ref {3.4} and Eq.(\ref {3.10''e1}), there exists an
algebra homomorphism  $\psi $ from ${\cal R} (G, g_i, \chi_i^{-1};
j\in J) $ to $diag ( kG[kQ_1^c, G,  g_i, \chi _i; i\in J])$ by
sending $x_j$ to $a_{g_j, 1}^{(j)}$ for any $j\in J$. By (ii),
$\psi$ is bijective. It is clear that $\psi$ is a graded braided
Hopf algebra isomorphism.

Let $R : = {\cal R}(G, g_i, \chi_i; i\in J)$.  Obviously, $R_{(1)}
\subseteq P(R)$. It is sufficient to show that any non-zero
homogeneous element $z\in R$, whose  degree $deg (z)$ is not equal
to 1, is not a primitive element. Obviously, $z$ is not a primitive
element when $deg (z)=0$. Now  $deg (z)>1$. We can assume, without
lost generality,  that there exist $\nu _1, \nu _2, \cdots, \nu _t
\in J $ such that
  $z = \sum _{\mid {\bf i} \mid = n} k_{ \bf i}
x^{\bf i}$, where $k_i \in k $, $x ^{\bf i} = x_{\nu_1}
^{i_1}x_{\nu_2} ^{i_2}\cdots x_{\nu_t} ^{i_t}$ with $ i_1 + i_2
+\cdots + i_t =n$.   It is clear
\begin {eqnarray} \label
{e611}  \Delta (z) = \sum _{\mid {\bf i} \mid =n} k_{\bf i}\Delta
(x^{{\bf i}}) = z\otimes 1 + 1\otimes z + \sum _{\mid {\bf i} \mid
=n} \ \sum_{ 0 \le {\bf j} \le {\bf i} , \  0 \ne {\bf j }\ne {\bf
i}} k _{\bf i} c_{\bf {i, j}}x^{\bf j}\otimes x^{\bf  i- j}.\end
{eqnarray} If $z$ is a primitive element, then   $\sum _{\mid {\bf
i} \mid =n} \ \sum_{ 0 \le {\bf j} \le {\bf i} , \  0 \ne {\bf j
}\ne {\bf i}} k _{\bf i} c_{\bf {i, j}}x^{\bf j}\otimes x^{\bf  i-
j}=0$. Since $c_{{\bf i}, {\bf j}} \not=0$, we have $k_{\bf i} =0$
for any ${\bf i}$ with $\mid {\bf i} \mid  =n$, hence $z=0$. We get
a contradiction. Thus $z$ is not a primitive element. This show
$R_{(1)} = P(R)$ and $R$ is a Nichols algebra.

(iv) and (v). Let $A:= {\cal T} (G, g_i, \chi_i ^{-1}; j\in J)$ and
$R = diag ( kQ^s( G, g_i, \chi _i; i\in J))$. Let $\psi$ be an
algebra homomorphism from ${\cal T} (G, g_i, \chi_i ^{-1}; j\in J)$
to $ diag ( kQ^s ( G,  g_i, \chi _i; i\in J))$ by sending $x_j$ to
$a_{g_j, 1}^{(j)}$.

It is clear that ${\cal T} (G, g_i, \chi_i^{-1}; j\in J)$ is a
$kG$-module algebra. Define a linear map $\nu $ from  $A \# kG$ to
$kQ^s$ by sending $x_j \# g$ to $a _{g_j, 1}^{(j)} \cdot g = \alpha
^+ (a _{g_j, 1}^{(j)} \otimes g)$ for any $g\in G, j \in J$. That
is, $\nu$ is the composition of
$$A \# kG \stackrel {\psi \otimes id} {\rightarrow } R \# kG
\stackrel {\alpha ^+} {\cong} kQ^s,
$$ where $\alpha ^+ = \mu _{kQ^s}( id \otimes \iota _0)$ (see
Theorem \ref {6}).  Define  a linear  map  $\lambda$ from $kG$ to
$A\#kG$ by sending $g$ to $1\# g$ for any $g \in G$ and another
linear map $\gamma$ from $kQ_1^c$ to $A\#kG$ by sending $a_{g_ih,
h}^{(i)}$ to $\chi _i^{-1}(h)x_i \# h $ for any $h \in G, i \in J.$
It is clear that $\gamma$ is a $kG$-bimodule homomorphism from
$kQ_1^c$ to $ _\lambda( A\#kG)_ \lambda$. Considering $kQ^s = T_{kG}
(kQ_1^c)$ and universal property of tensor algebra over $kG$, we
have that there exists an algebra homomorphism $\phi = T_{kG}
(\lambda, \gamma) $ from $kQ^s$ to $A\# kG$. Obviously, $\phi$ is
the inverse of $\nu$.
  Thus $\phi$ is bijective. It is easy to
check that $\phi$ is graded  Hopf algebra isomorphism. Obviously,
$\{  x_{\nu _1}
 x_{\nu _2} \cdots  x_{\nu _t} \# g \mid  \nu _j \in J,  j= 1, 2, \cdots
, t;  t \in {\mathbb Z}^+ \cup \{0\}, g\in G\}$ is a basis of ${\cal
T} (G, g_i, \chi_i^{-1}; j\in J) \# kG$. See
\begin {eqnarray*} && \nu ( x_{\nu _1}
 x_{\nu _2} \cdots  x_{\nu _t} \# g )\\
 &=&E_{\nu _1}\otimes _{kG}
 E_{\nu _2}\otimes _{kG}\cdots \otimes _{kG} E_{\nu _t} \cdot g \\
 &= & \chi _{\nu _1 } (g) \chi _{\nu _2 } (g) \cdots \chi _{\nu _t } (g) g\cdot E_{\nu _1}\otimes
_{kG}
 E_{\nu _2}\otimes _{kG}\cdots \otimes _{kG} E_{\nu _t}.
\end {eqnarray*}
Thus $ \{  g\cdot E_{\nu _1}\otimes _{kG}
 E_{\nu _2}\otimes _{kG}\cdots \otimes _{kG} E_{\nu _t} \mid  \nu _j \in J,  j= 1, 2, \cdots
, t;  t \in {\mathbb Z}^+ \cup \{0\}, g\in G\}$ is a basis of
$kQ^s$. It is easy to check that $\psi$ is graded braided Hopf
algebra isomorphism.

 (i) Considering (ii),  Lemma \ref {3.2} and Theorem \ref {6},
we complete the proof. \ \ $\Box$

\subsection {A characterization of multiple Taft algebras
}\label {s3.3}

In this subsection we characterize multiple Taft algebras by means
of elements in themselves.
\begin{Definition}\label{3.6} For a quantum weakly commutative   ${\rm ESC} (G, g_i, \chi _i; i\in
J)$,  let $A$ be the Hopf algebra to satisfy the following
conditions: {\rm(i)} $G$ is  a subgroup of $G(A)$;
 {\rm(ii)} there exists a linearly independent  subset $\{X_i \mid i\in J\}$ of $A$ such
 that $A$ is generated by set $\{X_i \mid i\in J\} \cup G$ as
 algebras;
 {\rm(iii)}  $X_j$ is $(1, g_j)$-primitive, i.e.,
$\Delta(X_j)=X_j\otimes 1+g_j\otimes X_j$, for any $j\in J;$
{\rm(iv)} $X_jg=\chi_j(g)gX_j$, for any $j\in J,$, $g\in G$;
 {\rm(v)}
$X_jX_i = \chi_j(g_i) X_iX_j$, for $i, j \in J$ with $i\not=j;$
 {\rm(vi)} $A_{(0)} \cap A_{(1)} =0$, where $A_{(0)}:=kG$ and
 $A_{(1)}$ is the vector space spanned by set $\{hX_i \mid  i \in J; h \in
 G\}$.
Furthermore, let  $J(A)$ denote   the ideal of $A$ generated by the
set
$$\{X_i^{N_i}  \mid N_i < \infty , i \in J \} $$ and
$H(G,g_i, \chi _i; i\in J )$  the quotient algebra $A/J(A)$.
\end{Definition}

\begin {Lemma} \label {3.5'} (See \cite [Lemma 3.3]{AS98b}  ) Let  $H$ be a Hopf algebra  with a linearly independent
set $ \{ x_{1} \dots, x_{t} \}$ and $G$ a subgroup of   $ G(H)$.
Assume that $g_{i} \in Z(G)$ and $\chi_{i} \in \hat G$ such that
$\Delta (x_i) = x_i \otimes 1+ g_i \otimes x_i$, \
 $x_ih=\chi_i(h)hx_i$,  for $i = 1, 2, \cdots,  t$, $h\in G.$ If the
intersection of $kG$ and span $\{hx_i \mid 1\le i \le t, h\in G\}$
is zero, then
\begin{eqnarray*}\{hx_1^{m_1} x_2^{m_2}\cdots x_t^{m_t} \mid 0 \leq
m_j<N_j, 1\le j \le t; h \in G \}.
\end{eqnarray*} is linearly independent, where $N_i$ is the order of
$q_i := \chi _i (g_i)$  \ ( $N_i = \infty $ when $q_i$ is not a root
of unit, or $q_i =1$ ) for $1 \le i \le t.$
\end {Lemma}

\begin{Proposition}\label{3.7} If    ${\rm ESC} (G, g_i, \chi _i; i\in
J)$ is  quantum weakly commutative, then
 $H(G;$ \ $ g_i, \chi_i; $ $ i\in J)$ and multiple Taft algebra $kG[kQ_1^c, G,  g_i,
\chi_i; $ $i\in J]$ are isomorphic as graded Hopf algebras.
\end{Proposition}
 {\bf Remark:} $H(G,g_i, \chi _i; i\in J )$ just is
$H(C, n, c, c^*, 0, 0)$ in \cite [ Definition 5.6.8 and Definition
5.6.15] {DNR01} with $G = C$, $J = \{1, 2, \cdots , t\}$, $1<n_i =
N_i < \infty$, $g_i =c_i, c_i^* = \chi _i$ for $i = 1, 2, \cdots,
t.$

\subsection {The relation between semi-path Hopf algebras and quantum enveloping
algebras}\label {3.4} If $0\not= q \in k$ and $0 \le i \le n < ord
(q) $ (the order of $q$), we set
$$ \left [ \begin {array} {c} n\\
i
\end {array} \right ]_q
 =  \frac{[n]_{q}!}{[i]_{q}![n - i]_{q}!}, \quad \hbox {where
}[n]_{q}! = \prod_{1 \le i \le n} [i]_{q}, \quad [n]_{q} =
\frac{q^{n} - q^{-n}}{q - q^{-1}}.$$

Let $B$ be a Hopf algebra and $R$  a braided Hopf algebra in ${}_B^B
{\cal YD}$. For convenience,  we denote   $r\# 1$ by $r$  and $1 \#
b$ by  $b$ in    biproduct $R \# B$ for any $r\in R,\  b \in B$.

\begin {Lemma} \label {3.7'}  Let $B$ be a Hopf algebra and $R$  a braided Hopf algebra in
${}_B^B {\cal YD}$ with $x_1, x_2 \in P(R)$ and $g_1, g_2 \in
Z(G(B))$. Assume that there exist $\chi_1, \chi _2 \in Alg(B,k)$
with $\sqrt{\chi _i (g_j)}\in k $ such that
 $$\delta _R(x_{j}) = g_{j}\otimes x_{j}, \ h\cdot  x_{j} = \chi_{j}(h ^{-1})x_{j},
$$
  for all $ h\in B, j =1, 2. $

(i) If $r$ is a positive integer and
  \begin {eqnarray} \label {quantume1'}
\begin {array} {cc}
 \chi _2 (g_1) \chi_1(g_2) \chi
_1 (g_1) ^{r-1} =1, \  \chi _1 (g_1) ^ {\frac {1}{2} (r-1)} \chi _2
(g_1) =1,\ r-1 < ord (\chi _1(g_1)),
\end {array}
\end {eqnarray} then
\begin {eqnarray}\label {quantume1'''}
 \sum_{m=0}^{r}(-1)^m
\left[\begin{array}{c}
r\\
m\\
\end{array}\right]_{\sqrt{\chi _1(g_1^{-1}})}
x_1^{r-m}x_2x_1^m \end {eqnarray} = $(ad_c x_1)^{r}x_2$ is a
primitive element of $ R$ and a $(1, g_1^rg_2)$-primitive element of
biproduct $R\#B$,  where $(ad_cx_1)x_2 = x_1x_2 - \chi _2
(g_1^{-1})x_2 x_1$.

  (ii) If $ \ \sqrt{\chi _1 (g_2)}  \sqrt{\chi _2(g_1)}=1$ and  $x _i g_j = \chi _i (g_j)g_jx_i$ for $i, j =1, 2$, then
\begin {eqnarray*}  \sqrt{\chi _2 (g_1)} x_1x_2- \sqrt{\chi
_1 (g_2)}x_2 x_1
\end {eqnarray*} is a primitive element of
$ R$  and
  \begin {eqnarray*}  \sqrt{\chi _2 (g_1)}x_1x_2- \sqrt{\chi
_1 (g_2)}x_2 x_1- \beta (g_1 g_2-1)
\end {eqnarray*} is a $(1, g_1g_2)$-primitive element of biproduct
$ R\#B$ for any $\beta \in k$.

\end {Lemma}

For an ${\rm ESC}(G, \chi_i, g_i; i\in J)$, we give the follows
notations:

(FL1) $N$  is a set and $J = \cup _{s\in N} (J_s \cup J_{s}')$ is a
disjoint union.

(FL2) There exists a bijection $\sigma : J^{(1)}\rightarrow J^{(2)}$
such that $\sigma \mid _{J_u}$ is a bijection from $J_u$ to $J_u'$
for $u \in N,$  where  $J^{(1)} := \cup _{u\in N} J_u $ and $J^{(2)}
:= \cup _{u\in N} J_u'.$

(FL3) There exists a  $J^{(1)}\times J^{(1)}$-matrix $A = (a_{ij})$
with $a_{ii} =2$ and non-positive integer $a_{ij} $ for any $i, j
\in J^{(1)}$ and $i \not= j$. For any $u \in N$, there exists an
integer $d_{i}^{(u)}$ such that $d_i ^{(u)}a_{ij} = d_j
^{(u)}a_{ji}$ for any $i, j \in J_u$.

(FL4) For any $u \in N$, there exists $0 \not= q_u \in k$  such that
$\chi _i (g_j) = q_u ^{-2d_i^{(u)}a_{ij}}$, $\chi _{\sigma (i)}
(g_j) = \chi _i ^{-1}(g_j)$  and $g_{\sigma (j)} = g_j $ for $i, j
\in J_u$.

(FL5) There exists $\xi _i\in G$ such that $\chi _{\sigma (i)} (\xi
_j) = \chi_i ^{-1}(\xi _j)$, $\xi _{\sigma (i)} := \xi _i ^{-1}$ and
$g_i= g_{\sigma (i)} := \xi _i ^2$  for any $i, j \in J_u, \ u \in
N;$ there exists a positive integer $r_{ij}$ such that  $r_{\sigma
(i), \sigma (j)}= r_{ij}$ for any $i, j \in J^{(1)} $ with $i\not= j
$.

(FL6) $
 \chi _j (g_i) \chi_i(g_j) \chi
_i (g_i) ^{r_{ij}-1} =1,$ $ \chi _i(\xi_j) = \chi _j(\xi _i), $ $
\chi _i (g_i) ^ {\frac {1}{2}  (r_{ij}-1)} \chi _j (g_i) =1,$ for
any  $i, j \in J_u $, $i \not= j,$ $ u\in N.$

(FL7) $G$ is a free commutative  group generated by  generator set
$\{\xi_i \mid i \in J^{(1)}\}$.

 An ${\rm ESC}(G, g_i, \chi _i; i\in J)$ is said to be a  local
FL-matrix type   (see \cite [P.4]{AS00}) if (FL1)--(FL4) hold. An
${\rm ESC}(G, g_i, \chi _i; i\in J)$ is said to be a  local FL-type
if (FL1), (FL2), (FL5) and (FL6) hold. An ${\rm ESC}(G, g_i, \chi
_i; i\in J)$ is said to be a  local FL-free type    if (FL1), (FL2),
(FL5), (FL6) and (FL7) hold. An ${\rm ESC}(G, g_i, \chi _i; i\in J)$
is said to be a  local FL-quantum group type    if (FL1)-- (FL4) and
(FL7) hold. If $N $ only contains one element  and $J^{(1)} = \{1,
2, \cdots, n\}$, then  we delete `local' in  the terms above.

Let ${\rm ESC}(G, \chi_i, g_i; i\in J)$ be a local FL-free type.
 Let $I$ be the ideal of
$kQ^s(G, g_i,\chi_i;i\in J)$ generated by the following elements:
\begin {eqnarray}\label {quantume5}
\chi _j (\xi_i)E_iE_j-\chi _i ^{-1}(\xi_j)E_j E_i- \delta_{\sigma
(i), j}\ \frac{g_i^2-1}{ \chi _i (\xi_i) -\chi _i (\xi _i) ^{-1}} ,
\ \hbox {for }          i \in  J_u, j \in J_u', u\in N;
\end {eqnarray}
\begin {eqnarray}\label {quantume6}
\sum_{m=0}^{r_{ij}}(-1)^m \left[\begin{array}{c}
r_{ij}\\
m\\
\end{array}\right]_{\chi _i(\xi_i^{-1})}
E_i^{r_{ij}-m}E_jE_i^m, \end {eqnarray} $  \hbox {  for any } i, j
\in J_u \hbox { or } i, j \in J_{u}',   i \not= j \hbox { and }
r_{ij}-1 < ord (\chi _i(g_i)) ; u\in N.$

  Let  $U$  be the  algebra generated  by set $\{K_i,   X_i \mid  i \in
J\}$ with relations \begin {eqnarray}\label
{quantume2}\begin{array}{c} X_iX_j-X_j X_i= \delta_{ \sigma(i), j} \
\frac{K _i^2- K_i ^{-2}}{\chi
_i(\xi_i) -\chi _i (\xi_i) ^{-1}},   \hbox { \ for any  } i \in J_u, j\in J_u'; u\in N ;\\
 K_iK_{\sigma (i)}=K_{\sigma (i)}K_i=1
 \hbox { \ for any  } i \in J ^{(1)} ;\\
\chi _j (\xi_i)K_iX_j=X_jK_i, \ \ K_iK_j-K_jK_i=0, \  \hbox { for any  } i, j  \in J ;\\
 \sum_{m=0}^{r_{ij}}(-1)^m \left[\begin{array}{c}
r_{ij}\\
m\\
\end{array}\right]_{\chi _i(\xi_i^{-1})}
X_i^{r_{ij}-m}X_jX_i^m =0,  \\
\hbox { for } i, j \in J_{u}  \hbox { or } i, j \in J_{u}',   i
\not= j \hbox { and }  r_{ij}-1 < ord (\chi _i(g_i)) ; u\in N.
\end {array} \end {eqnarray}
The comultiplication, counit and antipode of $U$ are defined by
\begin {eqnarray}\label {e3.811}\begin{array}{lll}
\Delta(X_j)=X_j\otimes K_{\sigma (j)}+K_j\otimes X_j,&
S(X_j)=- \chi _j (\xi _j)X_j,&\epsilon(X_i)=0,\\
\Delta(K_i)=K_i\otimes K_i,& S(K_j)=K_{\sigma (j)},
&\varepsilon(K_i)=1, \\
\Delta(X_{\sigma (j)})=X_{\sigma (j)}\otimes K_{\sigma
(j)}+K_j\otimes X_{\sigma (j)}, & S(X_{\sigma (j)})=- \chi _{\sigma
(j)} (\xi _j)X_{\sigma (j)}
\end{array}\end {eqnarray}
for any $j \in J^{(1)}, i \in J .$

 In fact,  $K_{\sigma (j)}= K_j ^{-1}$ in
$U$ for any $j \in J^{(1)}.$

\begin{Theorem}\label{8} Under notation above, if ${\rm ESC}(G, g_i, \chi _i; i\in
J)$ is  a  local FL-free type, then  $kQ^s(G, g_i,\chi_i;i\in J)/ I
\cong U$ as Hopf algebras.

\end {Theorem}

{\bf Proof}. For any $i, j \in J_u, u \in N, $  see $\chi _{\sigma
(i)} (\xi _{\sigma (j)}) = \chi _i (\xi _j) = \chi _j (\xi_i) = \chi
_{\sigma (j)} (\xi _{\sigma (i)}),$ $\chi _{\sigma (i)}(\xi_j) =
\chi _i ^{-1}(\xi _j) = \chi _j ^{-1} (\xi_i) = \chi _j (\xi
_{\sigma (i)}),$ $\chi _{i} (\xi_ {\sigma (j)}) = \chi _i (\xi _j)
^{-1} = \chi _j (\xi_i) ^{-1} = \chi _{\sigma (j)} (\xi _i ).$
Therefore,
\begin {eqnarray}\label {quantume13}
\chi _i (\xi _j) = \chi _j (\xi _i)
\end {eqnarray} For any $i, j \in J_u\cup J_u', u \in N.$
Obviously, for $i, j \in J_u'$, $i \not= j,$ $u\in N$, (FL6) holds.

We show this theorem by following several steps.

(i) There is a  algebra homomorphism $\Phi$ from $kQ^s$ to $U$ such
that $\Phi(\xi _i)= K_i$, $\Phi( a_{hg_i, h}^{(i)})= \Phi (h)K_iX_i$
and $\Phi( a_{hg_{\sigma (i)}, h}^{(\sigma (i))})= \Phi
(h)K_iX_{\sigma (i)}$ for all $h \in G$ and $i \in J ^{(1)}$.
Indeed, define algebra homomorphism $\phi: kG\rightarrow U$ given by
$\phi(\xi_i)=K_i$ for $i\in J$ and a $k$-linear map $\psi:
kQ_1^c\rightarrow U$ by $\psi(a_{gg_j,g}^{(j)})=\phi(g)K_jX_j $ and
$\psi(a_{gg_{\sigma (j)},g}^{(\sigma (j))})=\phi(g)K_jX_ {\sigma
(j)} $  for any $j \in J^{(1)},$ $ g \in G.$ For any $g, h\in G,
j\in J ^{(1)}$, see
$$\begin{array}{rcl}
\psi(h\cdot
a^{(j)}_{gg_j,g})&=&\psi(a^{(j)}_{hgg_j,hg})=\phi(hg)K_jX_j
=\phi(h)\phi(g)K_jX_j=\phi(h)\psi(a^{(j)}_{gg_j,g}),\\
\end{array}$$
and
$$\begin{array}{rcl}
\psi(a^{(j)}_{gg_j,g}\cdot h)&=&\chi_j(h)\psi(a^{(j)}_{hgg_j,hg})=\chi_j(h)\phi(hg)K_jX_j\\
&=&\phi(g)K_jX_j\phi(h) \ \ (\hbox { since } X_j \phi (h) = \phi (h) \chi _j (h) X_j)\\
&=&\psi(a^{(j)}_{gg_j,g})\phi(h).
\end{array}$$
Similarly, $\psi(h\cdot a^{(\sigma (j))}_{gg_{\sigma (j)},g}) =
\phi(h)\psi(a^{(\sigma (j))}_{gg_{\sigma (j)},g}) $ and $
\psi(a^{(\sigma (j))}_{gg_{\sigma (j)},g}\cdot h)= \psi(a^{(\sigma
(j))}_{gg_{\sigma (j)},g})\phi(h).$
 This implies that $\psi$ is a  $kG$-bimodule map from
$(kQ_1^c,g_i,\chi_i;i\in J)$ to $_{\phi}U_{\phi}$. Using the
universal property of tensor algebra over $kG$, we complete the
proof.

(ii) $\Phi (I)=0$. For any $i \in J_u, j \in J_{u}'$, see that
$$\begin{array}{cl}
&\Phi(\chi _j (\xi_i)E_iE_j-\chi _i (\xi_j)^{-1}E_j E_i-
\delta_{\sigma (i), j} \ \frac{g_i^2-1}{ \chi _i (\xi _i) -\chi _i
(\xi _i)
^{-1}})\\
 =& \chi _j (\xi_i) K_i X_iK_j^{-1} X_j-\chi _i (\xi_j)^{-1}
 K_j^{-1} X_jK_i X_i-
\delta_{\sigma (i),
j} \ \frac{K_i^4-1}{ \chi _i (\xi) -\chi _i (\xi _i) ^{-1}}\\
=& K_iK_j^{-1}X_iX_j-K_j^{-1}K_iX_jX_i- \delta_{\sigma (i),
j}\ \frac{K_i^4-1}{ \chi _i (\xi) -\chi _i (\xi _i) ^{-1}}\\
\\=&K_iK_j^{-1}(X_iX_j-X_jX_i-\delta_{\sigma (i),
j}
  \ \frac {K_i^3K_j-K_i^{-1}K_j}{ \chi _i (\xi _i) -\chi _i (\xi _i)^{-1} }) \\
=&K_iK_j^{-1}(X_iX_j-X_jX_i- \delta_{\sigma (i), j} \
\frac{K_i^2-K_i^{-2}}{ \chi _i (\xi _i) -\chi _i (\xi _i)^{-1} })=0.\\
\end{array}$$

For $i, j\in J_u,$ $i\not= j,$ see that
$$\begin{array}{cl}
&\Phi(E_i^{r_{ij}-m}E_jE_i^m)\\
=&(K_iX_i)^{r_{ij}-m}K_jX_j(K_iX_i)^m\\
=& \chi _i (\xi _i) ^{ \frac {1} {2} ( (r_{ij} -m) (r_{ij} -m -1) +
m (m-1) + 2 (r_{ij} -m) m)} \chi _i (\xi _j) ^{r_{ij} -m} \chi _j
(\xi _i)
^m  X_i^{r_{ij}-m}X_jX_i^m\\
=& \chi _i (\xi _i) ^{\frac {1}{2} (r_{ij}-1)r_{ij}} \chi _i (\xi
_j) ^{r_{ij}}X_i^{r_{ij}-m}X_jX_i^m \ \  ((\hbox {by }
(FL6))\\
\end{array}$$
and
$$\Phi(\sum_{m=0}^{r_{ij}}(-1)^m \left[\begin{array}{c}
r_{ij}\\
m\\
\end{array}\right]_{\chi _i(\xi_i^{-1})}
E_i^{r_{ij}-m}E_jE_i^m)=0.$$ Similarly, the equation above holds for
$i, j \in J_u'$, $i \not= j,$ $u\in N.$

By (ii), there exists an algebra homomorphism $\bar \Phi$ from
$kQ^s/I$ to $U$ such that $\bar \Phi (x +I) = \Phi (x)$ for any
$x\in kQ^s.$ For convenience, we will still use $x$ to denote $x
+I$.

(iii) It follows from the definition of $U$ that there exists a
unique algebra map $\Psi: U\rightarrow kQ^s/I$ such that
$\Psi(K_i)=\xi_i$, $\Psi(X_i)=\xi_i^{-1}E_i$ and $\Psi(X_{\sigma
(i)})=\xi_i^{-1}E_{\sigma (i)}$  for all $i \in J^{(1)}$. It is easy
to see that $\Psi\overline{\Phi}= id$ and $\bar \Phi \Psi =id$.

(iv) We  show that $I$ is a Hopf ideal of $kQ^s.$ It follows from
Theorem \ref {7} that $R:=diag (kQ^s)$ is a braided Hopf algebra
with $\delta ^- (E_i) =g_i \otimes E_i $, $h \rhd E_i = \chi _i
^{-1} (h) E_i$ for any $h \in G$, $i \in J.$ By Theorem \ref {6}, $R
\#B\cong kQ^s$ as Hopf algebras. It follows from Lemma \ref {3.7'}
that $I$ is a Hopf ideal of $kQ^s$.

 Obviously, $\Psi $ preserve the comultiplication and counit.
Since $kQ^s/I$ is a Hopf algebra, then  $\Psi$ is a Hopf algebra
isomorphism. \ \ $\Box$

\begin{Corollary}\label{3.14} The quantum enveloping algebra of a complex  semisimple Lie
algebra is isomorphic to a quotient of a semi-path Hopf algebra as
Hopf algebras.

\end {Corollary}
{\bf Proof.} Let  $k$ be the complex field and $L$ a complex
semisimple Lie algebra determined by $A = (a_{ij})_{n\times n}$. So
$A$ is a symmetrizable Cartan matrix with $d_i \in \{1, 2, 3\}$ such
that $d_ia_{ij} = d_j a_{ji}$ for any $i, j \in J ^{(1)} = \{1, 2,
\cdots, n\}$. Let $\mid \! N \! \mid =1, $  $J^{(2)} = \{n+1, n+2,
\cdots, n+n\}$, $J = J^{(1)}\cup J^{(2)}$ and $\sigma : J^{(1)}
\rightarrow J^{(2)}$ by sending $i$ to $i+n$. Let $G$ be a free
commutative group generated by  generator set $\{\xi_i \mid i \in
J^{(1)}\}$. Set $\xi _{\sigma (i)} = \xi _i^{-1}$,  $g_i= g_{\sigma
(i)} := \xi _i ^2$, $r_{ij}= r_{\sigma (i), \sigma (j)} = 1-a_{ij}$
for any $i, j \in J^{(1)}$, $i \not= j$. Define $\chi _i (\xi _j) =
q ^{-d_ia_{ij}}$ and $\chi _{\sigma (i)} (\xi _j) = \chi _i
^{-1}(\xi_j)$ for any $i, j \in J^{(1)}$, where $q$ is not a root of
1 with $0\not=q\in k$. It is easy to check (FL1)-- (FL7) hold, i.e.
${\rm ESC} (G, g_i, \chi_i; i\in J)$ is an FL-quantum group type.
 $U$ in Theorem \ref {8}
exactly  is the quantum enveloping algebra $U_q(L)$ of $L$ (see
\cite[p.218]{Mo93} or  \cite {Lu93}).  Therefore the conclusion
follows from Theorem \ref {8}. $\Box$

In fact, for any a generalized Cartan matrix $A$, we can obtain an
FL-quantum group type ${\rm ESC} (G, g_i, \chi_i; i\in J)$ as in the
proof above. By the way, if ${\rm ESC} (G, g_i, \chi_i; i\in J)$ is
a Local FL-type, then there exist the  Hopf ideals, which generated
by (\ref {quantume5}) and (\ref {quantume6}),  in co-path Hopf
algebra $kQ^c (G, g_i, \chi_i; i\in J)$ and multiple Taft algebra
$kQ^c (G, g_i, \chi_i; i\in J)$, respectively.

\section {Appendix}\label {s4}

\begin{Proposition}\label {1.12} If $(kQ_1^c, G, r, \overrightarrow \chi, u )$ is a $kG$-Hopf bimodule, then
  $(kG)^*$-coactions on the $(kG)^*$-Hopf bimodule $(kQ_1^a,
G, r, \overrightarrow \chi, u)$ are given by
$$\delta^-(a^{(i)}_{y,x})=\sum _{h \in G} p _h \otimes  a
^{(i)} _{h^{-1}y, h^{-1}x},\ \ \delta^+(a^{(i)}_{y,x})=\sum_{h\in G}
\chi _C^{(i)}(\zeta_{\theta}(h^{-1})^{-1})a^{(i)}_{yh^{-1},
xh^{-1}}\otimes
 p_h$$
where $x, y\in G$ with $x^{-1}y=g^{-1}_{\theta}u(C)g_{\theta}$,
$\zeta_{\theta}$ is given by {\rm(\ref{e0.3})}, $C\in{\mathcal
K}_r(G)$, $i\in I_C(r)$ and $p_h\in(kG)^*$ is defined by
$p_h(g)=\delta_{h,g}$ for all $g, h\in G$.
\end{Proposition}

Now, we give  an interesting quantum combinatoric formula by means
of multiple Taft algebras. Let $n$ be a positive integer and $S_n$
be the symmetric group on the set $\{1,2,\cdots,n\}$. For any
permutation $\sigma\in S_n$, let $\tau(\sigma)$ denote the number of
reverse order of $\sigma $, i.e., $\tau(\sigma )=|\{(i,j)\mid 1\leq
i<j\leq n, \sigma(i)>\sigma(j)\}|$. For any $0\not= q\in k$, let
$S_n(q):=\sum_{\sigma\in S_n}q^{\tau(\sigma)}$.

\begin{Lemma}\label{3.10} In $kQ^c(G, r, \overrightarrow \chi, u)$, assume ${g} \in
{\mathcal K}_r (G)$ and  $j\in I_{\{g\}}(r)$. Let $q:= \chi
^{(j)}_{\{g\}} (g)$.
 If $i_1, i_2, \cdots, i_m$ be non-negative integers, then
$$\begin{array}{rcl}
a^{(j)}_{g^{i_m +1},g^{i_m }}\cdot a^{(j)}_{g^{i_{m-1}
+1},g^{i_{m-1}}}\cdot\cdots\cdot a^{(j)}_{g^{i_1 +1},g^{i_1 }}
&=&q^{\beta_m+\frac{m(m-1)}{2}}S_m(q^{-1})P^{(j)}_{g^{\alpha_m}}(g,m)
\end{array}$$
where $\alpha _m = i_1 + i_2 + \cdots + i_m $, $P^{(j)}_h(g,m) =$ \
$ a^{(j)}_{g^mh,g^{m-1}h}a^{(j)}_{g^{m-1}h, g^{m-2}h}\cdots
a^{(j)}_{gh,h}$, $\beta_1=0$ and $\beta_m=\sum_{j
=1}^{m-1}(i_1+i_2+\cdots+i_j )$ if $m>1$.

\end{Lemma}
{\bf Proof.} Now let $a_l=a^{(j)}_{g^{i_l+1},g^{i_l}}$ for
$l=1,2,\cdots,m$ and consider the sequence of $m$ arrows $A=(a_m,
\cdots, a_1)$. We shall use the notations of \cite[p.247]{CR02}. For
any $\sigma\in S_m$ and $1\leq l\leq m$, let $m(\sigma,l)=|\{i\mid
i<l,\sigma(i)<\sigma(l)\}|$. Then
$$m(\sigma,l)=(l-1)-|\{i\mid
i<l,\sigma(i)>\sigma(l)\}|,$$ and hence
$$\begin{array}{rcl}
\sum_{l=1}^mm(\sigma,l)&=&\frac{m(m-1)}{2}-\sum_{l=1}^m|\{i\mid
i<l,\sigma(i)>\sigma(l)\}|\\
&=&\frac{m(m-1)}{2}-\tau(\sigma).\\
\end{array}$$
Now it follows from \cite[Proposition 3.13]{CR02} that
$$\begin{array}{rcl}
a_m\cdot a_{m-1}\cdot\cdots\cdot a_1&=& \sum_{\sigma\in S_m}
A_m^\sigma\cdots  A_2^\sigma A_1^\sigma\\
&=&\sum_{\sigma \in S_m}(\prod_{l=1}^m q^{i_1+i_2+\cdots
+i_{\sigma(l)-1}+m(\sigma,l)})P^{(j)}_{g^{\alpha_m}}(g,m)\\
&=&q^{\beta_m+\frac{m(m-1)}{2}}S_m(q^{-1})
P^{(j)}_{g^{\alpha_m}}(g,m).\ \ \Box
\end{array}$$\\

 Let $G\cong{\bf  Z}$ be the infinite cyclic group with
generator $g$ and let $r$ be a ramification of $G$ given by
$r_{\{g\}}=1$ and $r_{\{g^n\}}=0$ if $n\not=1$. Let $0\not=q\in k$.
Define $\chi_{g} ^{(1)}\in\widehat G$ by $\chi_{g}^{(1)}(g)=q$. Then
$(G, r, \overrightarrow{\chi}, u)$ is an ${\rm RSC}$. Using  Lemma
\ref{3.4} and Lemma \ref {3.10}, one gets the following result.

\begin{Example}\label{3.11}
{\rm(i)}\ For any  \ $ 0\not=q\in k$,
$(m)_q!=q^{\frac{m(m-1)}{2}}S_m(q^{-1}) $, \ and \ $ S_m(q) =$ $
(m)_{\frac {1}{q}}!q^{\frac{m(m-1)}{2}}$, where $m$ is a positive
integer.

{\rm(ii)} \  Assume $q$ is a primitive $n$-th root of unity with
$n>1$. Then
$S_m(q)=0$ if $m\geq n$, and $S_m(q)\not=0$ if $0<m <n$.\\
\end{Example}

\vskip 0.5cm

\noindent{\large\bf Acknowledgement}:  We would like to thank
referees and Prof. N. Andruskiewitsch for comments and suggestions.
The first two authors were financially supported by the Australian
Research Council, and the author Chen was supported by NSF of China
(10471121) and Sino-German project (GZ310). S.C.Z thanks the
Department of Mathematics, University of Queensland for hospitality.

\end{document}